\documentclass{amsart}
\usepackage{amscd,amssymb}
\numberwithin{equation}{section}

\newtheorem {Theorem} 			{Theorem}
\newtheorem {Question}		{Question}
\newtheorem*{Question*}		{Question}
\newtheorem {varTheorem}                {Theorem}
\newenvironment {Theorem'}
        {\begin{varTheorem}{\hspace{-3.5mm}}{\bf '}{\hspace{3.5mm}}}
        {\end{varTheorem}}
\newtheorem {RefTheorem}[equation]     	{Theorem}  	
\newtheorem {Lemma}[equation]     	{Lemma}  	
\newtheorem {Proposition}[equation]	{Proposition}  

\newtheorem {Corollary}[equation]	{Corollary}

\newtheorem {corollary}[equation]	{Corollary}

\theoremstyle{definition}

\newtheorem {Remark}[equation]		{Remark}
\newtheorem {Example}[equation]		{Example}

\newcommand{\pr} {\smallskip\noindent{\bf Proof\,\,}}

\newcommand     {\printname}[1] {}

\newcommand{\labell}[1] {\label{#1}\printname{#1}}

\def    \SO	{\operatorname{SO}}
\def    \SU	{\operatorname{SU}}
\def    \Sq	{\operatorname{Sq}}

\def	\sign	{\operatorname{sign}}
\def	\mod	{\ {\operatorname{ mod }}\ }
\def	\Tilde	{\widetilde }
\def	\Hat	{\widehat}

\def	\to	{\longrightarrow}

\def	\C	{{\mathbb C}}
\def	\R	{{\mathbb R}}

\def	\Q	{{\mathbb Q}}
\def	\Z	{{\mathbb Z}}
\def	\CP	{{\mathbb C}{\mathbb P}}

\def	\ft	{{\mathfrak t}}
\def	\fg	{{\mathfrak g}}

\def	\SO	{\operatorname{SO}}

\def	\pr	{\operatorname{pr}}

\begin{document}

%%%%%%%%%%%%%%%%%%%%%%%%%%%%%%%%%%%%%%%%%%%%%%%%%%%%%%%%%%%%%%%%%%%%%%%%%%%%%%%
%%%%%%%%%%%%%%%%%%%%%%%%%%%%%%%%%%%%%%%%%%%%%%%%%%%%%%%%%%%%%%%%%%%%%%%%%%%%%%%

\title[On  a symplectic generalization of Petrie's conjecture]
{On  a symplectic generalization of Petrie's conjecture}

\author{Susan Tolman}

\address{Department of Mathematics,  University of Illinois at
Urbana-Champaign, 
Urbana, IL 61801}
\email{stolman@math.uiuc.edu}

\thanks{The author was partially supported by National Science Foundation
Grant DMS \#07-07122.}

\begin{abstract}
Motivated by the Petrie conjecture,  we consider the following
questions:
Let  a circle act in a Hamiltonian fashion on a compact symplectic manifold $(M,\omega)$
which satisfies $H^{2i}(M;\R) = H^{2i}(\CP^n,\R)$ for all $i$. 
Is $H^j(M;\Z) = H^j(\CP^n;\Z)$ for all $j$?  Is the total
Chern class of $M$ determined by the cohomology ring $H^*(M;\Z)$?
We answer these questions in the six dimensional case by showing that
$H^j(M;\Z)$ is equal to $H^j(\CP^3;\Z)$ for all $j$,
by proving that only four cohomology rings can arise, and by computing the total Chern class 
in each case.  
We also prove that there are no exotic actions.
More precisely, if $H^*(M;\Z)$ is isomorphic to $H^*(\CP^3;\Z)$
or $H^*(\Tilde{G}_2(\R^5);\Z)$, then 
the representations at the fixed components 
are compatible with one of the  standard actions;
in the remaining two case, the representation 
is strictly determined by the cohomology ring.
Finally,  our results suggest a natural question: do the
remaining two cohomology rings actually arise?
This question is closely related to some  interesting problems
in symplectic topology, such as embeddings of ellipsoids.
\end{abstract}

\maketitle

%\tableofcontents

%%%%%%%%%%%%%%%%%%%%%%%%%%%%%%%%%%%%%%%%%%%%%%%%%%%%%%%%%%%%%%%%%%%%%%%%%%%%%%
%%%%%%%%%%%%%%%%%%%%%%%%%%%%%%%%%%%%%%%%%%%%%%%%%%%%%%%%%%%%%%%%%%%%%%%%%%%%%%

\section{Introduction}

In the early 1970's, Ted Petrie
wanted to address two related 
fundamental  questions: Given a compact Lie group
$G$ and a  manifold $M$, does $M$ admit a  
$G$ action\footnote
{ We shall always assume that our actions are non-trivial.}?
If so, how many different actions can we find?

One of his important insights was that these questions are much
more tractable when $M$ is a
{\bf homotopy projective space}, that is,
a simply connected manifold so that $H^*(M;\Z) = H^*(\CP^n;\Z)$ as rings,
or equivalently a manifold which is homotopy equivalent to $\CP^n$.

For example, the first key step in answering these questions
is understanding the relationship
between  the tangent bundle
near the fixed components and the global invariants of $M$.
Petrie  proved that  if the circle  acts on a homotopy
projective space with isolated fixed points, then the Pontrjagin
classes are determined by the representations at the fixed points
\cite{Pe1}.

Motivated by  this and other evidence, he stated
what is now known as the {\bf Petrie conjecture}:
if a homotopy projective
space $M$ admits a circle action, then the Pontrjagin classes
of $M$ are {\bf standard}, that is, agree with the Pontrjagin
classes of $\CP^n$ itself.
Although this conjecture has not been  resolved in general,
it has motivated a good deal of research.
In particular, 
it has been proven if  $M$ is at most eight dimensional 
\, \cite{Dej, Ja},
if  $M$ admits an invariant almost complex structure  
whose first Chern class 
is at least $\frac{1}{2} \dim (M) + 1$ times the generator of $H^2(M;\Z)$ \cite{Ha},
and in many other special cases \cite{Des, Ma, Mu, Pe2, TsWa,
Wan, Yo}.

We are interested in addressing the 
analogous  questions for symplectic manifolds:
Given a compact Lie group
$G$ and a symplectic  manifold $(M,\omega)$, does $M$ admit a  Hamiltonian
$G$ action?  
If so, how many different actions can we find?

In the symplectic case,
several additional tools are available.
For example, there is
an almost complex structure $J \colon T(M) \to T(M)$ which is
{\bf compatible}  with $\omega$,  i.e.,
$\omega(J \cdot , \cdot)$ is a Reimannian metric.
Moreover, the set of such structures is contractible, and so 
there is a well-defined total Chern class  $c(M) \in H^*(M;\Z)$.
Additionally, the components of the moment map $\Phi \colon M \to \fg^*$ 
are Morse-Bott functions with extremely nice properties; see \S\ref{s:background}.

Therefore, instead of insisting  that our symplectic manifold
$(M,\omega)$ be a  homotopy
projective space,  we merely assume that 
$H^{2i}(M;\R) = H^{2i}(\CP^n;\R)$ for all $i$.
We prove that if the circle 
acts on such a manifold in a Hamiltonian fashion with isolated fixed points, 
then both the cohomology ring and total Chern class are determined
by the representations at the fixed points;
see Corollary~\ref{cor:isom} and Remark~\ref{rmk:cob}. 
This leads to the following questions.

\begin{Question} \rm
Consider a Hamiltonian  circle action  on a symplectic manifold
$(M,\omega)$ which satisfies  $H^{2i}(M;\R) = H^{2i}(\CP^n;\R)$ for all $i$.
Is $H^j(M;\Z) = H^j(\CP^n;\Z)$ for all $j$?
Is the total Chern class $c(M)$  completely determined by
the cohomology ring $H^*(M;\Z)$?
\end{Question}

Our first main theorem  answers this question affirmatively in
the $6$-dimensional case.
In fact, we are able to  show that only a few possible rings arise.

\begin{Theorem}\labell{thm:1}
Let the circle  act  on a $6$-dimensional
compact symplectic manifold $(M,\omega)$ with 
moment map $\Phi \colon M \to \R$.
If  $H^2(M,\R) = \R$, then one  of the following four statements is true:
\begin{itemize}
\item [(A)]
$H^*(M;\Z) = \Z[x] /(x^4) $
and  $c(M)= 
1 + 4x + 6x^2 + 4x^3.$
\item [(B)]
$H^*(M;\Z) = \Z[x,y] /(x^2 - 2 y, y^2)$
and  $c(M)= 
1 + 3x + 8y + 4xy.$
\item [(C)] $H^*(M;\Z) = \Z[x,y] /(  x^2 - 5 y, y^2)$
and  $c(M)= 
1 + 2x + 12y + 4xy.$
\item [(D)] $H^*(M;\Z) = \Z[x,y] /(  x^2 - 22 y, y^2)$
and  $c(M)= 
1 + x + 24y  + 4xy.$
\end{itemize}
In each case, $x$ has degree $2$ and $y$ has degree $4$.
\end{Theorem}

\begin{Remark}\labell{rmk:betti}
Since $H^2(M;\R) = \R$ and $M$ is a $6$-dimensional symplectic manifold,
Poincare duality implies that $b_0 = b_2 = b_4 = b_6 = 1$, where $b_i = \dim (H^i(M;\R))$ denotes
the $i$'th Betti number.
If the fixed set is discrete,  then this immediately implies that $H^j(M;\Z) = H^j(\CP^3;\Z)$
for all $j$, that is, that $b_1 = b_3 = b_5 = 0$ and the cohomology
is torsion free; see \S\ref{s:background}.
However, when the fixed set is not discrete, this fact is  somewhat surprising;
it follows from the analysis in \S\ref{s:notdiscrete}. 
\end{Remark}

%\begin{Remark}
%The conclusions of Theorems \ref{thm:1}, \ref{thm:2}, and \ref{thm:3} also hold
%if a circle acts on a $6$-dimensional compact symplectic manifold $(M,\omega)$ which satisfies
%$H^1(M;\R) = 0$ and $H^2(M;\R) = \R$,
%because in this case every symplectic action is Hamiltonian.
%\end{Remark}

\begin{Remark}
By Corollary~\ref{pi1} below, the  manifold $M$ described above is
simply connected.
If we assume that  statement (A) is true then $w_2(M) = 0$; 
therefore, Wall's theorem implies that 
$M$ and $\CP^3$ are diffeomorphic \cite{Wal}.
Similarly, any two manifolds which satisfy statement (C) 
must be diffeomorphic.
\end{Remark}

Petrie was able to construct  {\bf exotic} circle actions on 
projective spaces, that is, actions
so that the induced representations at the 
fixed points do not agree
with those of  any circle subgroup $S^1 \subset \SU(n+1)$.
Our second main theorem is 
that, in contrast,  in the $6$-dimensional symplectic case their are no
exotic actions. More precisely, the representations
at the fixed components are either 
strictly determined by the cohomology ring
or are compatible with one of the standard actions described 
below.

\begin{Example}\labell{ex:a}
Given $n > 1$, let $\CP^n$ denote the projective
space of lines in $\C^{n+1}$.
Since this $2n$ dimensional manifold naturally arises
as a coadjoint orbit of $\SU(n+1)$, it inherits a symplectic form  $\omega$
and a Hamiltonian $\SU(n+1)$ action.
Hence, every circle subgroup $S^1 \subset \SU(n+1)$ induces a
Hamiltonian circle action on $\CP^n$.
\end{Example}

\begin{Example}\labell{ex:b}
Given $n > 1$, let $\Tilde G_2(\R^{2n+1})$ denote the Grassmannian of oriented
$2$-planes in $\R^{2n+1}$. Since this $4n - 2$ dimensional manifold 
naturally arises
as a coadjoint orbit of $\SO(2n+1)$, it inherits a symplectic form  $\omega$
and a Hamiltonian $\SO(2n+1)$ action.
Hence, every circle subgroup $S^1 \subset \SO(2n+1)$ induces a
Hamiltonian circle action on $\Tilde{G}_2(\R^{2n+1})$.
\end{Example}

Given any subgroup  $H \subset S^1$,
let  $M^H$ denote the set of points fixed by $H$.
Each component $N \subset M^H$ is a symplectic manifold which inherits
a symplectic circle action with moment map $\Phi|_N$.
If $H \neq \{e\}$, we call each component $N$ of $M^H$ which is
{\em not } fixed by $S^1$ an {\bf isotropy submanifold}. 
Each two-dimensional isotropy submanifold is a sphere which contains
exactly two isolated fixed points;
we call these {\bf isotropy spheres}.
We can now state our second main theorem, which 
is an immediate consequence of 
Propositions
\ref{prop:notdiscrete}, 
\ref{prop:simple}, and \ref{prop:multiple},
and the remarks subsequent to each.

\begin{Theorem}\labell{thm:2}
Let the circle  act faithfully\footnote{
A group $G$ acts {\bf faithfully} on $M$ if for every non-trivial $g \in G$
there exists $m \in M$ so that $g \cdot m \neq m$.}
on a $6$-dimensional
compact symplectic manifold $(M,\omega)$ with 
moment map $\Phi \colon M \to \R$.
If  $H^2(M,\R) = \R$,
then one of the following four statements is true:\footnote{
Throughout this paper,  the symbol $\cong$ implies that
the two sides are  equivariantly isomorphic  complex vector bundles.}
\begin{itemize}
\item [(A)]
There is a subgroup $S^1 \subset \SU(4)$ and   an orientation preserving
diffeomorphism
$f \colon M^{S^1} \to \bigl(\CP^3 \bigr)^{S^1}$
so that
$T (M)|_{M^{S^1}} \cong f^*\left(
T \bigl(\CP^3 \bigr) \big|_{ ( \CP^3 )^{S^1} }
\right). $ 
\item [(B)]
There is a subgroup $S^1 \subset \SO(5)$ and an orientation preserving  diffeomorphism
$f \colon M^{S^1} \to \Tilde{G}_2(\R^5)^{S^1} $  
so that 
$T(M)|_{M^{S^1}} \cong 
f^*\left(T \bigl(\Tilde{G}_2(\R^5) \bigr) 
\big|_{\Tilde{G}_2(\R^5)^{S^1}}\right) $. 
\item [(C)]
The fixed set consists of four points; 
the weights at these points are
%%S $\{1,2,3\}$, $\{1,-1,4\}$, $\{1,-1,-4\}$, and $\{-1,-2,-3\}$.
$$\{1,2,3\}, \{1,-1,4\}, \{1,-1,-4\},\mbox{ and } \{-1,-2,-3\}.$$
\item [(D)]
The fixed set consists of four points; 
the weights at these points are
%%S $\{1,2,3\}$, $\{1,-1,5\}$, $\{1,-1,-5\}$, and $\{-1,-2,-3\}$.
$$\{1,2,3\}, \{1,-1,5\}, \{1,-1,-5\}, \mbox{ and } \{-1,-2,-3\}.$$
\end{itemize}
Moreover, $M$ contains a pair
of isotropy spheres which intersect in two points in cases (C) and (D), but
not in cases (A) or (B).
\end{Theorem}

\begin{Remark}
In case (A),
$M$ is cobordant to
$\CP^3$ (with some multiple of the standard symplectic form) as a stable-complex Hamiltonian $G$-space.  In fact,
we may assume that $f$ is a symplectomorphism;  see Remarks~\ref{rmk:cob} and \ref{rmk:cob2}.
Similar comments apply in each case.
\end{Remark}

We can convert any non-trivial circle action into an effective
circle action by quotienting out the subgroup which acts trivially.
Therefore, Theorem~\ref{thm:1} follows immediately
from Theorem~\ref{thm:2} above and  
Corollaries \ref{cor:isom} and \ref{cor:isol}; see
Example~\ref{ex:non}.

In each of the cases described above,
each component of the fixed set is simply 
connected. By \cite{Li1}, this implies that Theorem \ref{thm:2}
has the following Corollary.

\begin{Corollary}\labell{pi1}
Let the circle  act faithfully on a $6$-dimensional
compact symplectic manifold $(M,\omega)$ with 
moment map $\Phi \colon M \to \R$.
If  $H^2(M,\R) = \R$,
then $M$ is simply connected.
\end{Corollary}

Finally, we can describe
the {\bf equivariant cohomology} of $M$;
by definition, this is
$H^*_{S^1}(M) = H^*(M \times_{S^1} S^\infty).$
For example, if $p$ is a point then
$H^*_{S^1}(p;\Z) = H^*(\CP^\infty;\Z) = \Z[t].$
The projection map $\pi \colon M \times_{S^1} S^{\infty} \to \CP^\infty$ 
induces a pull-back map 
\begin{equation}\labell{pi}
\pi^* \colon H^*(\CP^\infty;\Z) \to H_{S^1}^*(M;\Z);
\end{equation}
thus, $H^*_{S^1}(M;\Z)$ is a $H^*(\CP^\infty;\Z)$ module.
Moreover, the inclusion
$\iota \colon M^{S^1} \to M$ induces a restriction map
$\iota^* \colon H^*_{S^1}(M;\Z) \to H^*_{S^1}(M^{S^1};\Z)$;
define $$H^*_{S^1}(M;\Z)|_{M^{S^1}} = \iota^* (H^*_{S^1}(M;\Z)).$$
Finally, let  $c^{S^1}(M) \in H_{S^1}^*(M;\Z)$ 
denote the total equivariant Chern class of $M$.

The theorem below follows  immediately from
Theorem~\ref{thm:2} and Corollaries~\ref{cor:eqisom} and \ref{cor:eqisol};
see Example~\ref{ex:eq}.
Note that, since $H^*(M^{S^1};\Z)$ has no torsion,
the image 
$H^*_{S^1}(M;\Z)|_{M^{S^1}}$
naturally determines 
the equivariant cohomology ring itself;
see \S\ref{s:background}.

%\begin{Theorem}
%Let the circle $S^1$ act faithfully on a $6$-dimensional
%compact symplectic manifold $(M,\omega)$ with moment map 
%$\Phi \colon M \to \R$.
%If  $H^2(M,\R) = \R$,
%then  one of the following four statements is true:
%\begin{itemize}
%\item [(A)] There is an isomorphism 
%$\hat{f} \colon H^*_{S^1}(M;\Z) \to H^*_{S^1}(\CP^3;\Z)$ so that
%$\hat{f}(c(M)) = c(\CP^3)$
%for some $S^1 \subset (S^1)^3$.
%\item 
%[(B)] There is an isomorphism 
%$\hat{f} \colon H^*_{S^1}(M;\Z) \to H^*_{S^1}(\Tilde{G}_2(\R^5);\Z)$ 
%so that 
%$\hat{f}(c(M)) =  c(\Tilde{G}_2(\R^5))$
%for some $S^1 \subset (S^1)^2$.
%\item  [(C)] 
%$ \Z[x,y,t]/ 
%(x^2 - 5 y + x t, y^2 + 2 x y t + 6 y t^2),$
%and
%$c(M) = 1  + 2 x + 6t + 12 y+ 12 x t + 11 t^2 +4 xy
%+ 32 y t + 10 x t^2 + 6 t^3.$
%\item [(D)]
%$H^*_{S^1}(M) = \Z[x,y,t]/ (x^2 - 22 y + x t, 
%y^2 +  xy t + 6 y t^2),$
%and
%$$c(M) = 1 + x + 6 t + 24 y + 12 x t + 11 t^2
%+ 4 xy + 68 y t + 11 x t^2 + 6 t^3.$$
%\end{itemize}
%In each case $x$ has degree $2$ and $y$ has degree $4$.
%\end{Theorem}

\begin{Theorem}\labell{thm:3}
Let the circle  act faithfully on a $6$-dimensional
compact symplectic manifold $(M,\omega)$ with moment map 
$\Phi \colon M \to \R$.
If  $H^2(M,\R) = \R$,
then  one of the following four statements is true:
\begin{itemize}
\item [(A)] 
There is a subgroup $S^1 \subset \SU(4)$ and diffeomorphism
$f \colon M^{S^1} \to \bigl(\CP^3 \bigr)^{S^1}$
so that
\begin{gather*}
H_{S^1}^*(M;\Z)\big|_{M^{S^1}}
= 
f^* \bigl(H_{S^1}^*(\CP^3 ;\Z)\big|_{(\CP^3)^{S^1}} \bigr)  
\mbox{ and }   \\
c^{S^1}(M)\big|_{M^{S^1}}
= 
f^*\bigl(c^{S^1}(\CP^3)\big|_{(\CP^3)^{S^1}}\bigr) 
.\end{gather*}
\item [(B)]
There is a subgroup $S^1 \subset \SO(5)$ and  diffeomorphism
$f \colon M^{S^1} \to \Tilde{G}_2(\R^5)^{S^1} $  
so that
\begin{gather*}
H_{S^1}^*(M;\Z) \big|_{M^{S^1}}
= 
f^*\left(H_{S^1}^*\bigl(\Tilde{G}_2(\R^5) ;\Z\bigr) 
\big|_{\Tilde{G}_2(\R^5)^{S^1}}\right)  
\mbox{ and } \\
c^{S^1}(M)|_{M^{S^1}}
= 
f^*\left(c^{S^1}\bigl(\Tilde{G}_2(\R^5)\bigr) \big|_{\Tilde{G}_2(\R^5)^{S^1}}\right) 
.\end{gather*}
\item[(C)] The fixed set consists of four points: $p_0, p_1, p_2$ and $p_3$.
As a $H^*(\CP^\infty;\Z) = \Z[t]$ module, $H_{S^1}^*(M;\Z)$ is generated by
$1, \alpha_1,\alpha_2$, and $\alpha_3$,
where  
\begin{gather*}
\alpha_1|_{p_1} = t, \ \alpha_1|_{p_2} = 5t,\  \alpha_1|_{p_3} = 6 t, \
\alpha_2|_{p_2} = 4 t^2, \\
 \alpha_2|_{p_3} = 6 t^2, \
 \alpha_3|_{p_3} = 6 t^3,  
\ \mbox{and} \ \alpha_i|_{p_j} = 0 \ \forall \ j < i;  \ \mbox{moreover},\\
c^{S^1}(M)|_{p_0} =  1 + 6t + 11 t^2 + 6 t^3, \  \
c^{S^1}(M)|_{p_1} =  1 + 4t  -  t^2 - 4 t^3,  \\
c^{S^1}(M)|_{p_2} =  1 - 4t  -  t^2 + 4 t^3 \ \ \mbox{and}  \ \ 
c^{S^1}(M)|_{p_3} =  1 - 6t + 11 t^2 - 6 t^3.
\end{gather*}
\item[(D)] The fixed set consists of four points: $p_0, p_1, p_2$ and $p_3$.
As a $H^*(\CP^\infty;\Z) = \Z[t]$ module, $H_{S^1}^*(M;\Z)$ is generated by
$1, \alpha_1,\alpha_2$, and $\alpha_3$,
where 
\begin{gather*}
\alpha_1|_{p_1} = t, \ \alpha_1|_{p_2} = 6t,\  \alpha_1|_{p_3} = 12 t, \
\alpha_2|_{p_2} = 5 t^2, \\
 \alpha_2|_{p_3} = 6 t^2,  \
 \alpha_3|_{p_3} = 6 t^3, \
\ \mbox{and} \ \alpha_i|_{p_j} = 0 \ \forall \ j < i; \ \mbox{moreover},  \\
c^{S^1}(M)|_{p_0} =  1 + 6t + 11 t^2 + 6 t^3, \ \ 
c^{S^1}(M)|_{p_1} =  1 + 5t  -  t^2 - 5 t^3,  \\
c^{S^1}(M)|_{p_2} =  1 - 5t  -  t^2 + 5 t^3, \ \  \mbox{and} \ \ 
c^{S^1}(M)|_{p_3} =  1 - 6t + 11 t^2 - 6 t^3, \\
\end{gather*}
\end{itemize}
\end{Theorem}

Note that, in the case that the action is semifree and  there is no four dimensional
fixed component, these three theorems are due to Li \cite{Li2}.  

\subsection*{Open questions}

In these theorems, 
the first two cases correspond to Examples~\ref{ex:a} and \ref{ex:b},
but the last two cases do not correspond to any known examples\footnote{
Since this paper was originally submitted, 
McDuff has used symplectic techniques to
construct  manifolds corresponding to the last two cases \cite{Mc2}.
In fact, as she points out, both manifolds (which are K\"ahler) 
were already known.}.
This raises the following natural question:

\begin{Question}\labell{Q2}\rm
Do there exist examples exhibiting properties (C) or (D) of Theorem
2?  
More precisely, let $l = 4$ or $5$.
Does there exist a Hamiltonian circle action on a compact symplectic (alternatively,
K\"ahler) manifold 
so that the fixed set consists of four points 
with weights $\{1,2,3\}$, $\{1,-1,l\}$, $\{1,-1,-l\}$,
and $\{-1,-2,-3\}$?
\end{Question}

We do not know the answer to this question.
These manifolds cannot be
ruled out by any of the techniques used in this paper;
see also Remarks~\ref{Wu} and ~\ref{HBJ}.

Moreover, this question seems to be related to
interesting problems  in symplectic topology.
To see this, we need to introduce some more notation.
Given a  $n +1$-tuple  of natural numbers $\mathbf {k} = (k_0, \ldots, k_n)$, 
consider the weighted projective space of type $\mathbf{k}$,
$$ \CP^n( \mathbf{k}) = 
S^{2n+1} /
{(z_0, \ldots, z_n) \sim (\lambda^{k_0} z_0, \ldots, \lambda^{k_n} z_n)}
;$$ 
let $\overline{\CP}^n(\mathbf{k})$ denote the same manifold with the
opposite orientation.
Let $\alpha_{\mathbf{k}} \in H^2 \bigl( \CP^n(\mathbf{k}) \bigr)$
denote the first Chern class of
the tautological  circle bundle
$S^{2n+1} \to  \CP^n(\mathbf{k})$.
Finally, given real number $a$ and $b$, 
define  the ellipsoid
$$ E(a,b) = \left\{ (x_1,x_2) \in \C^2  \left| \
\tfrac{1}{a} x_1^2 + \tfrac{1}{b} x_2^2 \leq  1 \right. \right\}  .$$

Now  suppose that a manifold $(M,\omega)$ satisfying the conditions
of Question 2 does exist.
By  Corollary~\ref{cor:eqisom} 
(see Example~\ref{ex:eq}) and Lemma~\ref{le:extend}
 -- after possibly rescaling $\omega$ --
there exists a moment map $\Phi \colon M \to \R$ so that 
$$\Phi(p_0) = -6, \ 
\Phi(p_1) = - l, \ \Phi(p_2) = l , \ \mbox{and} \ \Phi(p_3) = 6,$$
where $p_i$ is the unique fixed point of index $2i$ for all $i$ such
that $0 \leq 2i \leq 6$.
By \cite{Go},  this implies that
for all $\kappa \in (-l,l)$, the reduced space
$M_\kappa = \Phi^{-1}(\kappa)/S^1$ is diffeomorphic to 
to the connected sum 
$$X = \CP^2(1,2,3) \# \overline{\CP}^2(1,1,l).$$
Moreover, let $\omega_\kappa \in \Omega^2(M_\kappa)$ denote
the reduced symplectic form;
the cohomology class $[\omega_\kappa] \in H^2(M_\kappa)$ 
is the unique class
so that
$[\omega_\kappa]\big|_{\CP^1(2,3)} = (6 + \kappa) \alpha_{(2,3)}$
and
$[\omega_\kappa]\big|_{\overline{\CP}^1(1,l)} = - (l + \kappa) \alpha_{(1,l)}$.
Here, the inclusions of $\CP^1(2,3)$ and $\CP^1(1,l)$ into
$X = M_\kappa$ are induced by the natural inclusions
$\CP^1(2,3) \subset \CP^2(1,2,3)$
and $\overline{\CP}^1(1,l) \subset \overline{\CP}^2(1,1,l).$
In particular, an affirmative answer to Question 2 implies an
affirmative answer to the  question below.
%\begin{quote}

\begin{Question} \rm
Given  any $\lambda < 2$,
is there a symplectic (K\"ahler) form 
\begin{gather*}
\omega_\lambda \in 
\Omega^2\bigl( \CP^2(1,2,3) \# \overline{\CP}^2(1,1,l) \bigr)
\quad \mbox{such that} \\
[\omega_\lambda]\big|_{\CP^1(2,3)} = (6 + l) \alpha_{(2,3)}
\ \mbox{and} \
[\omega_\lambda]\big|_{\overline{\CP}^1(1,l)} = - \lambda l \alpha_{(1,l)}?
\end{gather*}
\end{Question}

%for all $\kappa \in (-6,-l)$ the reduced space
%$M_\kappa$ is diffeomorphic to  $\CP^2(1,2,3)$
%and $[\omega_\kappa] =  (\kappa + 6) \Tilde{c_1}(\mathcal{O}_{\mathbf{n}}(1))$.

As in the manifold case, we can  construct such a symplectic form
if we can find the symplectic embeddings  described below.
\begin{Question} \rm
Given any $\lambda < 2$,
is there  a symplectic embedding 
$$E(\lambda, \lambda l)  \hookrightarrow
E\bigl(\tfrac{6+l}{3}, \tfrac{6+l}{2}\bigr)?$$
\end{Question}
Unfortunately, although symplectic embeddings have been extensively studied,
this particular question does not seem to follow easily from known results
\cite{S}\footnote{
Again, McDuff has published new results on this question
since this paper was originally submitted \cite{Mc1}; she uses these
in \cite{Mc2}.}.
In particular, it cannot be ruled out by volume constraints;
$\frac{\text{volume}( E(2, 2 l))}  
{\text{volume}\bigl( E\bigl(\frac{6+l}{3}, \frac{6+l}{2}\bigr)\bigr)}$ is 
equal to $\frac{24}{25}$ if $l = 4$ and to $\frac{120}{121}$
if $l = 5$.

%\begin{Remark}
%Let the circle act on a compact symplectic manifold $(M,\omega)$
%so that statement (A) of Theorem 2 holds.
%Rescale the symplectic form so that the cohomology class
%$[\omega]$ generated $H^2(M;\Z)$.
%
%If there is a $4$-dimensional fixed component, 
%$M$ is equivariantly symplectomorphic to $\CP^3$ with its standard
%symplectic form; see \S7.
%
%Otherwise, we cannot prove this claim.
%However,
%$M^{Z_k}$ and  $(\CP^{3})^{\Z_k}$ are equivariantly symplectomorphic 
%for all $k > 1$.
%To see this, note that,
%after translating if necessary,
%$F$ and $f(F)$ have the same moment image for every
%fixed component $F$.
%Moreover, $F$ and $f(F)$ have the same area
%for  every $2$-dimensional fixed component $F$.
%Finally, by Theorem~\ref{thm:2},
%$M$ does not contain two isotropy spheres which intersect in two points.
%It is straightforward to check that this implies that
%two components $F$ and $F'$ lie in the same component
%of $M^{\Z_k}$ exactly if $f(F)$ and $f(F')$
%lie in the same component of $(\CP^3)^{\Z_k}$.
%Since each component of $M^{\Z_k}$ is at most $4$-dimensional,
%the claim follows immediately from Yael's classification
%of Hamiltonian circle actions on $4$-dimensional compact symplectic manifolds 
%\cite{Kar}.

%A similar remark applies if statement (B) of Theorem~\ref{thm:2} is true
%(except that the fixed set never has a $4$-dimensional component).
%\end{Remark}

We conclude this section with a brief overview of this paper.
In \S\ref{s:background}, we introduce some background material and
establish our notation.
In \S\ref{s:arbitrary}, we prove a few useful results which hold
in arbitrary dimensions. As a consequence, we prove that
Theorem 1 and Theorem 3 follow immediately from Theorem 2.
In \S\ref{s:notdiscrete}, we return to the six-dimensional case,
proving Theorem 2 in the case that the fixed set is not discrete.
We spend the remainder of the paper proving this theorem  in the case that
the fixed set {\em is} discrete. To do so, we first define
a labeled multigraph associated to $M$ in \S\ref{s:graph}, 
and then prove Theorem 2 in the cases  that the associated multigraph
is simple and not simple in \S\ref{s:simple} and \S\ref{s:notsimple},
respectively.

\subsection*{Acknowledgments}
I would like to thank Jonathan Weitsman  for inspiring these
results by introducing me to the Petrie conjecture.
I would also like to thank Yael Karshon and Dusa McDuff for useful discussions,
and the referee for suggesting many
improvements to the exposition.

\section{Background}\labell{s:background}
%\comment{Check on the first time symbols are used.}
%\comment{Need something like -- If $\beta$ is of degree $i$
%and vanishes all all fixed points of index less that $i$, then
%it must be the one of degree $i$.}
%\comment{Maybe put claim about  $\omega + t \Phi$ being an equivariant form
%in this section?}

In this section, we introduce some background material
and establish our notation.

Let $M$ be a compact manifold.
A {\bf symplectic} form on $M$ is a closed, non-degenerate
two-form $\omega \in \Omega^2(M)$.
A circle action on $M$ is {\bf symplectic} if it preserves $\omega$.
A symplectic circle action is {\bf Hamiltonian} if there exists a {\bf moment map},
that is,  a map $\Phi \colon M \to \R$ such that
$$- d \Phi = \iota_{\xi_M} \omega,$$
where $\xi_M$ is the vector field on $M$ induced by the circle action.
Since $\iota_{\xi_M} \omega$ is closed, every symplectic action is
Hamiltonian if $H^1(M;\R) = 0$.

Let the circle  act on a compact symplectic manifold $(M,\omega)$
with moment map $\Phi \colon M \to \R$.
Since the set of compatible almost complex structures
$J \colon T(M) \to T(M)$ is contractible,  there is a well-defined
multiset of integers, called {\bf weights}, associated to  each fixed point $p$.
Indeed, for any fixed component $F$,
the tangent bundle $T(M)|_F$ naturally splits into subbundles -- one
corresponding to each weight.

The moment map  $\Phi \colon M \to \R$
is a Morse-Bott function  whose 
critical set  is the fixed point set.
Moreover, the negative normal bundle at $F$
is the sum of the subbundles of $T(M)|_F$
with negative weights.  
In particular, the index of a fixed component $F$ is $2 \lambda_F$,  
where $\lambda_F$ is the number of negative weights in $T_p M$ for any $p \in F$
(counted with multiplicity). 
More interestingly, let $e_{S^1}(N^-_F) \in H^{2 \lambda_F}_{S^1}(F)$ denote the equivariant
Euler class of the negative normal bundle  at $F$.
If $p \in M^{S^1}$ is an isolated point,
then $e_{S^1}(N^-_p) = \Lambda_p^- t^{\lambda_p}$, where
$\Lambda_p^- \in \Z \smallsetminus \{0\}$ is the product of the negative (integer) weights at $p$.
More generally, for any fixed component $F$,
we can naturally identify $H^*_{S^1}(F) = H^*(F \times \CP^\infty)$ with $H^*(F)[t]$.
Under this identification, $e_{S^1}(N^-_F)$ is a polynomial
in $t$; the highest degree term is $ \Lambda_F^- t^{\lambda_F}$.
Therefore,
as Atiyah and Bott pointed out,
$e_{S^1}(N^-_F)$ is not a zero divisor in $H^{2 \lambda_F}_{S^1}(F;\Q)$ for  any  fixed component $F$.

Kirwan  uses this idea to prove three remarkable theorems:
``perfection'', ``injectivity'', and ``formality'' 
\cite{Ki}. 
Let $R = \Z$ if the fixed set is torsion free,
that is, $H^*(M^{S^1};\Z)$ has no torsion; otherwise, let $R = \Q$.
(See \cite{TWquot} for comments on the integral case.)
Let $F$ be any fixed component, and
let $M^\pm = \Phi^{-1}(-\infty, \Phi(F) \pm \epsilon)$,
where $\epsilon > 0$ is sufficiently small.
Since 
$e_{S^1}(N^-_F) \in H^*_{S^1}(F;R)$ is not a zero divisor,
the natural restriction $H_{S^1}^*(M^+,M^-;R) \to H_{S^1}^*(F;R)$ is an injection.
Therefore, the long exact sequence in equivariant cohomology for the pair
$(M^+,M^-)$ breaks into short exact sequences
\begin{equation}
\label{short}
 0 \to H_{S^1}^j(M^+,M^-;R) \to H_{S^1}^j(M^+;R)
\to H_{S^1}^j(M^-;R) \to 0.
\end{equation}
If $H_{S^1}^j(M^-;R)$ is a free group, this implies immediately that
$$H_{S^1}^j(M^+;R) =  H_{S^1}^j(M^+,M^-;R) \oplus H_{S^1}^j(M^-;R).$$
By induction,  $H_{S^1}^j(M;\R)$ is a free group and
the moment map is an equivariantly  perfect Morse-Bott function; in fact,
$$ H_{S^1}^j(M;R) = \bigoplus_{F \subset M^{S^1}} H_{S^1}^{j - 2\lambda_F}(F;R),$$
where the sum is over all fixed components.
Similarly, by  induction and \eqref{short}, the  restriction
map $\iota^* \colon H^*_{S^1}(M;R) \to H^*_{S^1}(M^{S^1};R)$ is an injection.
Hence, every equivariant cohomology class is  determined by
its restriction to the fixed point set.
Finally, restriction induces a natural map of exact sequences
\begin{equation*}
\minCDarrowwidth20pt
\begin{CD}
0  @>>> H^j_{S^1}(M^+,M^-;R) @>>> H^j_{S^1}(M^+;R)   @>>> H^j_{S^1}(M^-;R)  @>>> 0  \\
@. @VVV @VVV @VVV @. \\
\dots  @>>> H^j(M^+,M^-;R) @>>> H^j(M^+;R)   @>>> H^j(M^-;R)  @>>> \dots.  \\
\end {CD} 
\end{equation*}
Moreover, the restriction map from $H^*_{S^1}(M^+,M^-;R)$ to $H^*(M^+,M^-;R)$
is surjective because $H^*_{S^1}(F;R) = H^*(F;R)[t]$.
Hence, by an easy diagram chase, if the restriction map
from $H^*_{S^1}(M^-;R)$ to $H^*(M^-;R)$ is surjective,  then 
so is the restriction map from $H^*_{S^1}(M^+;R)$ to $H^*(M^+;R)$;
moreover, the long exact sequence in cohomology for the
pair $(M^+,M^-)$ breaks into short exact sequences
\begin{equation*}
 0 \to H^j(M^+,M^-;R) \to H^j(M^+;R)
\to H^j(M^-;R) \to 0.
\end{equation*}
By induction, $\Phi$ is a perfect Morse-Bott function  the restriction map
$H^*_{S^1}(M;R) \to H^*(M;R)$ is a surjection.
By Leray-Hirsch, this implies that 
the kernel of this map is the ideal generated by $\pi^*(t)$,
where $t \in H^2(\CP^\infty;R)$ is the generator.  (See \eqref{pi}.)
Hence, if  want to compute the (ordinary) cohomology of $M$,
it is enough to determine the equivariant cohomology of $M$
as a $H^*(\CP^\infty;R)$ module.
Nearly identical arguments prove the following closely related 
proposition.

\begin{Proposition}\labell{basis}
Let the circle action on a compact symplectic manifold $(M,\omega)$
with moment map $\Phi \colon M \to \R$.
If the fixed set is torsion-free, let $R = \Z$; otherwise let
$R = \Q$.

Given a fixed component $F$ and a class $u \in H^i(F;R)$,
there exists a class $\alpha \in H^{i + 2\lambda_F}_{S^1}(M;R)$ so that 
\begin{enumerate}
\item $\alpha|_F =  u \,  e_{S^1}(N^-_F)$, and
\item $\alpha|_{F'} = 0$ for all other fixed components 
$F'$ with $\Phi(F' ) \leq \Phi(F)$.
\end{enumerate}
%
%Moreover, consider $\alpha_1,\dots,\alpha_N \in H_{S^1}^*(M;R)$,
%Assume that for each $j$,  $\alpha_j \in H^{i_j}_{S^2}(M;R)$
%satisfies (1) and (2) above for some
%fixed component $F_j$ and class $u_j \in H^{i_j + 2 \lambda_{F_j}}(F_j;R)$.
%If $\{u_j\}$ is a basis for $H^*(M;R)$, then
%$\{\alpha_j\}$ is a basis for $H_{S^1}^*(M;R)$ 
%as a $H^*(\CP^\infty;R)$ module. \\
%
Moreover, let $\{u_j\}$ be a basis for $H^*(M^{S^1};R)$, where each $u_j \in H^{i_j}(F_j;R)$ for
some fixed component $F_j$.  If $\alpha_j \in H^{i_j + 2 \lambda_{F_j}}_{S^1}(M;R)$ and $u_j$
satisfy  (1) and (2) above for each $j$, then 
$\{\alpha_j\}$ is a basis for $H_{S^1}^*(M;R)$ 
as a $H^*(\CP^\infty;R)$ module. 
\end{Proposition}

\begin{corollary}\labell{cor:basis}
Let the circle action on a compact symplectic manifold $(M,\omega)$
with moment map $\Phi \colon M \to \R$.
If the fixed set is torsion-free, let $R = \Z$; otherwise let $R = \Q$.
Consider $\beta \in H_{S^1}^*(M;R)$ and $c \in \R$ so that  $\beta|_{F'} = 0$ 
for all fixed components $F'$ such that $\Phi(F') < c$.
\begin{enumerate}
\item If $c = \Phi(F)$ for some fixed component $F$,
then $\beta|_F$ is a multiple of $e_{S^1}(N^-_F)$.
\item More generally,
let $\{u_j\}$ be a basis for $H^*(M^{S^1};R)$, where  each $u_j \in H^{i_j}(F_j;R)$ for
some fixed component $F_j$. 
Assume that $\alpha_j \in H^{i_j + 2 \lambda_{F_j}}_{S^1}(M;R)$ and $u_j$
satisfy  (1) and (2) above for each $j$.  
Then $$\beta = \sum_{\Phi(F_j) \geq c} x_j\, \alpha_j,$$
where $x_j \in H^*(\CP^\infty;R)$ for all $j$.
Here, the sum is over all $j$ such that $\Phi(F_j) \geq c$. 
\end{enumerate}
\end{corollary}

\begin{proof}
By the proposition above, we can write $\beta = \sum_j x_j \alpha_j$, where
here the sum is over all $j$.  If $x_j = 0$ for all $j$ such that $\Phi(F_j)
< c$, then the second claim holds.  Moreover, if $\Phi(F) = c$ for some
fixed component $F$, then properties (1) and (2) 
together imply that $\beta|_F$ is a multiple of $e_{S^1}(N^-_F)$. 

Otherwise, there exists $j$ so that $\Phi(F_j)  < c$ and $x_j \neq 0$
but $x_k = 0$ for  $k$ such that $\Phi(F_k) < \Phi(F_j)$.
By properties (1) and (2), this implies that $\beta|_{F_j} \neq 0$.
Since $\Phi(F_j) < c$, this contradicts the assumption.
\end{proof}

%We can also use Morse theory to prove several more facts which are 
%peripherally relevant to out story:  A symplectic circle action is Hamiltonian
%exactly if there is a fixed component $F$ of index $0$.
%Moreover, if the action is Hamiltonian, then
%there is
%a unique minimal fixed component $F$ and maximal fixed component $F'$.

%We can also use Morse theory to prove another fact: 
%Given a Hamiltonian circle action on a symplectic manifold,
%there is
%a unique minimal fixed component $F$ and maximal fixed component $F'$.

The projection $\pi \colon M \times_{S^1} S^\infty \to \CP^\infty$
induces a natural push-forward map $\pi_* \colon H_{S^1}^*(M;\Z)
\to H^*(\CP^\infty;\Z)$.  Since this map is given by ``integration over
the fiber,'' we will usually denote it by the symbol $\int_M$.
We will need the following theorem, due to Atiyah-Bott and
Berline-Vergne \cite{AB,BV}.

\begin{RefTheorem} \labell{ABBV}
Let the circle act a compact manifold $M$.
Fix $\alpha \in H_{S^1}^*(M;\Q)$. As elements of $\Q(t)$,
$$\int_M \alpha = \sum_{F \subset M^{S^1}} \int_F \frac{ \alpha|_F}{e_{S^1}(N_F)},$$
where the sum is over all fixed components,
and $e_{S^1}(N_F)$ denotes the equivariant Euler class of the normal
bundle to $F$.
\end{RefTheorem}

\begin{Remark}\labell{rmk:ABBV}
If $p \in M^{S^1}$ is an isolated fixed point, and the (integer) weights at $p$
are $\xi_1,\ldots,\xi_n$ (repeated  with multiplicity), then
$c^{S^1}_i(M)|_p = \sigma_i(\xi_1,\ldots,\xi_n) t^i$
where $\sigma_i$ is the $i$'th elementary symmetric polynomial
and $t$ is the generator of $H^2_{S^1}(p;\Z)$.
For example, $c^{S^1}_1(M) = \sum \xi_i t$ 
and  $e_{S^1}(N_p) = c^{S^1}_n(M)|_p =  \left( \prod \xi_j \right) t^n.$
Hence,
$$\int_p \frac{ c^{S^1}_i|_p}{e_{S^1}(N_p)}= \frac{ \sigma_i(\xi_1,\ldots,\xi_n) }
{  \prod \xi_j } t^{i -n}.$$

If $\Sigma$ is a fixed
surface  of genus $g_\Sigma$ instead, then -- 
since every vector bundle over a surface splits --
the normal bundle to $\Sigma$ is the direct sum of line bundles
with equivariant Chern classes $\xi_1 t + a_1 u , \dots, \xi_{n-1} t + a_{n-1} u,$
where $\xi_1,\dots,\xi_{n-1}$ are the non-zero weights at $\Sigma$ (repeated
with multiplicity), 
$u$ is the positive generator of $H^2(\Sigma;\Z)$,
and the $a_i$'s are integers.
Since $e_{S^1}(N_\Sigma) = 
\left( \prod_j \xi_j t \right) \left( 1  + \frac{u}{t} \sum_j \frac{a_j  }
{\xi_j }  \right)$
and $1|_\Sigma = 1$, 
$$
\int_\Sigma \frac{ 1|_\Sigma}{e_{S^1}(N_\Sigma)} =
\int_\Sigma \frac{1}
{ \left( \prod_j \xi_j t \right) \left( 1  +  \frac{u}{t} \sum_j \frac{a_j}{\xi_j}  \right) }
=  
\frac 
{ \int_\Sigma \left( 1  - \frac{u}{t} \sum_j \frac{a_j}{\xi_j} \right)   }
{\prod_j \xi_j t  } 
= 
- \frac{  \sum_j \frac{a_j}{\xi_j}} 
{ \prod_j \xi_j  } 
t^{-n}.
$$
Similarly, since
$c_1^{S^1}(M)|_\Sigma = 
c_1(\Sigma) + u \sum_j  a_j  + t \sum_j \xi_j  $,
\begin{align*}
\int_\Sigma \frac{ c_1^{S^1}(M)|_\Sigma}{e_{S^1}(N_\Sigma)} &= 
\frac
{
\int_\Sigma \left( c_1(\Sigma) + u \sum_j  a_j  + t \sum_j \xi_j  \right)
\left(  1  - \frac{u}{t} \sum_j \frac{a_j}{\xi_j}   \right) 
}
{ \prod_j \xi_j t  }  
\\
&= 
\frac{ \int_\Sigma
\left( c_1(\Sigma)  -  u \sum_{i \neq j}  \frac{a_j \xi_i }{\xi_j} 
 + t \sum_j \xi_j  \right)}
{ \prod_j \xi_j t  }  
\\
&= 
\frac
{2 (1 -  g_\Sigma)  - \sum_{i \neq j} \frac{a_j \xi_i }{\xi_j} }
{ \prod_j \xi_j} t^{1-n}.
\end{align*}
Finally since , $(c_1^{S^1}(M))^2 - 2 c_2^{S^1}(M) = t^2 \sum_j \xi_j^2  + 
2  u t   \sum_j \xi_j a_j $,

\begin{align*}
\int_\Sigma \frac{ (c_1^{S^1}(M))^2|_\Sigma -  2 c_2^{S^1}(M)|_\Sigma}{e_{S^1}(N_\Sigma)} & = 
\frac
{ \int_\Sigma
\left( t^2  \sum_j \xi_j^2  + 2  u t \sum_j \xi_j a_j  \right)
\left(  1  - \frac{u}{t} \sum_j \frac{a_j}{\xi_j}   \right) 
}
{ \prod_j \xi_j t  }  
\\
&= 
 \frac
{ \int_\Sigma \left(
 u t  \sum_j a_j \xi_j  -  u t  \sum_{i \neq j} \frac{a_j \xi_i^2}{\xi_j}   + t^2 \sum_j  \xi^2_j  \right) }
{\prod_j \xi_j t  }  
\\
&=
\frac
{  \sum_j  a_j \xi_j - \sum_{i \neq j} \frac{a_j \xi_i^2 }{\xi_j} }
{ \prod_j \xi_j}  t^{2-n}.
\end{align*}
\end{Remark}

%\begin{example}
%For example, consider the circle action on $\CP^2$ given by
%$\lambda \cdot [x_0,x_1,x 2] = [x_0,x_1, \lambda x_2]$.
%There are two fixed components; the point $p = [0,0,1]$ and
%the sphere $\Sigma = [0,x_1,x_2]$. The equivariant Euler classes
%of the normal bundles are $e(N_p) = t^2$ and $ e(N_\Sigma) = t + u$,
%where $u$ is the positive generator of $H^2(\Sigma;\Z)$.
%Then 
%$$\int_{\CP^2} 1 =  \int_p \frac{1}{t^2} + \int_\Sigma \frac{1} {t + u}
%= \frac{1}{t^2} + \frac{1}{t} \int_\Sigma  \frac{1}{1 + \frac{u}{t}}
%= \frac{1}{t^2} + \frac{1}{t} \int_\Sigma  \left( 1 - \frac{u}{t} + \frac{u^2}{t^2} +  \dots  \right)
%= \frac{1}{t^2} - \frac{1}{t^2} = 0,$$
%as required.
%\end{example}

Finally, we will need the following very simple lemmas.

\begin{Lemma}\labell{modulo}
Let the circle act on a compact symplectic manifold $(M,\omega)$.
Let $p$ and $p'$ be fixed points which lie in the same component
$N$ of $M^{\Z_k}$, for some $k > 1$.
Then the  $S^1$-weights at $p$ and at $p'$ are equal modulo $k$.
\end{Lemma}

\begin{proof}
Since $\Z_k$ fixes $N$, the weights of the representation of
$\Z_k$ on the tangent space $T_q M$ are the same for all $q \in N$.
Moreover, if $q \in N$ is fixed by the circle action, then the
weights for the $\Z_k$ action on $T_q M$ are exactly the reduction
modulo $k$ of the weights for the circle action.
\end{proof}

\begin{Lemma}\labell{le:extend}
Let the circle act on a compact symplectic manifold $(M,\omega)$
with moment map $\Phi \colon M \to \R$. There exists an
equivariant extension $w \in H^2_{S^1}(M;\R)$ of $[\omega] \in H^2(M;\R)$
so that $w|_{F'} = [\omega|_{F'}] - \Phi(F') t$ for all fixed components $F'$.
\end{Lemma}

\begin{proof}
Take $w = [\omega - \Phi t]$ in the Borel model for equivariant cohomology.
\end{proof}
                                                                               
\begin{Example}\labell{proj} \rm
Let the torus $(S^1)^3 \subset \SU(4)$  act on $\CP^3$ by
$$ (\lambda_1,\lambda_2,\lambda_3) \cdot [x_0,x_1,x_2,x_3] =
[x_0, \lambda_1 x_1,\lambda_2 x_2,\lambda_3 x_3].$$
Let $e_1,e_2$ and $e_3$ denote the standard basis for
the weight lattice $(\Z^3)^*$.
The fixed points are $[1,0,0,0]$, $[0,1,0,0]$, $[0,0,1,0]$, and 
$[0,0,0,1]$;
the weights at these points are 
$\{e_1,e_2,e_3\}$, $\{-e_1, e_2 - e_1, e_3 - e_1\}$,
$\{-e_2, e_1 - e_2 , e_3 - e_2\}$,
and 
$\{-e_3 ,e_1 - e_3, e_2 - e_3 \}$, respectively.
\end{Example}

\begin{Example}\labell{gras}\rm
Let the torus $(S^1)^2 \subset \SO(5)$ act on  $\Tilde{G}_2(\R^5)$
induced by the  $S^1$ action on $\R^5 = \R \times \C^2$
given by
$$(\lambda_1,\lambda_2)  \cdot (t,x_1,x_2) = 
(t,  \lambda_1 x_1, \lambda_2 x_2  ).$$
Let $e_1$ and $e_2$  denote the standard basis for
the weight lattice $(\Z^2)^*$.
The fixed points are the planes 
$\{ (y_1,\ldots,y_5) \in \R^5 \mid y_1 = y_4 = y_5 = 0 \}$
and $\{ (y_1,\ldots,y_5) \in \R^5 \mid y_1 = y_2= y_3 = 0 \},$
with either orientation;
the weights  at these points  are
$\{e_1,e_1 + e_2, e_1 - e_2\}$,
$\{-e_1,-e_1 + e_2, -e_1 - e_2\}$,
$\{e_2,e_2 + e_1, e_2 - e_1\}$, and
$\{-e_2,-e_2 + e_1, -e_2 - e_1\}$.
\end{Example}

\begin{Remark}\labell{Wu}
The cases (C) and (D) described in Theorem~\ref{thm:1} are  consistent with
Wu's theorem, so we cannot use this theorem to rule out such manifolds.
Wu's theorem states that total Steifel-Whitney class
$w(M)$ of a connected manifold $M$ is equal to $\Sq(v(M))$, where
$\Sq \colon H^*(M;\Z_2) \to H^*(M;\Z_2)$ is the Steenrod square operator
and $v(M) \in H^*(M;\Z_2)$ is the unique class
so that $v(M) \cup x = \Sq(x)$ for all  $x \in H^*(M;\Z_2)$.  (See \S18.8 in \cite{Hu}.)
%Similarly, if $M$ is oriented and $p$ is prime,
%then the image of the total Pontrjagin class $p(M)$
%under the coefficient homomorphism $H^*(M;\Z) \to H^*(M;\Z_p)$ 
%is $\P(v(M))$,  
%where  $\P \colon H^*(M;\Z_p) \to H^*(M;\Z_p)$ is
%the Steenrod power operator
%and now $v(M) \in H^*(M;\Z_p)$
%is the unique class so that $v(M) \cup x = \P(x)$ 
%for all  $x \in H^*(M;\Z_p)$.
For example, if
$H^*(M;\Z_2) = \Z_2 [x]/ (x^4)$, as in case (C),
then $\Sq^2(x^2) = 0$ and so Wu's theorem implies that $w(M) = 0$.  
Since $w(M)$ is the image of $c(M)$ under the
coefficient homomorphism $H^*(M;\Z) \to H^*(M;\Z_2)$, 
this is satisfied both cases.
Similarly, if $H^*(M;\Z_2) = \Z_2 [x,y]/ (x^2,y^2)$,
as in case (D),
then either  $\Sq^2(y) = 0$ and $w(M) = 1$,
or $\Sq^2(y) =  xy$  and $w(M) = 1 + x$.
The latter statement is satisfied in both cases.
%Finally, if $H^*(M;\Z_3) = \Z_3 [x]/ (x^4)$, as in both cases,
%then $\P^1(x) = x^3$ and  $\P(v(M)) = 1 + x^2$.
%Since $p_1(M) = c_1(M)^2 -  2c_2(M)$,
%this  also holds in both cases.
\end{Remark} 

\begin{Remark}\labell{HBJ}
Similarly, cases (C) and (D) are consistent with the fixed point
formula for the Hirzebruch genus.
Let $M$ be a compact almost complex manifold; let
$\chi_y$ denote the Hirzebruch genus corresponding to the power
series 
$$\frac{ x \left(1 + y e^{-x(1+y)} \right)} {1 - e^{-x(1+y)}}.$$
On the one hand, if $M$ is $6$-dimensional then a direct
calculation shows that
$$\chi_y(M) = 
\frac{1}{24} (1 + y - y^2 - y^3) \int_M   c_1(M) c_2(M)
+ \frac{1}{2}( -y + y^2) \int_M c_1(M).$$
On the other hand, if a circle acts on $M$ then by \cite{HBJ},
$\chi_y(M)  = \sum_{F} (-y)^{\lambda_F} \chi_y(F),$
where the sum is of all fixed components.
In particular,  
in cases (C) and (D), 
$$\chi_y(M) = 1 - y + y^2 - y^3.$$
Since 
$\int_M c_1(M) c_2(M) = 24$ and $\int_M c_3(M) = 4$
in both cases, these formulas agree.
\end{Remark}

\section{Arbitrary dimensions}
\labell{s:arbitrary}

We begin by exploring some of the consequences of our central
assumptions in {\em arbitrary} dimensions.
More precisely, let 
a circle act in a Hamiltonian fashion  on compact symplectic manifold  $(M,\omega)$;
assume that $H^i(M;\Z) = H^i(\CP^n;\Z)$ for all $i$.
Our  main result is that the fixed point data determines the
(equivariant) cohomology ring and Chern classes; see Proposition~\ref{prop:isom}.
We also show, in Proposition~\ref{cor:order}, that the index of the fixed components
determines the order of their moment images.

Roughly speaking, our first result states that critical 
components near the minimum must have low index. 

\begin{Lemma}\label{le:order}
Let the circle act on a compact symplectic manifold $(M,\omega)$
with moment map $\Phi \colon M \to \R$.
Given any fixed component $F$, 
$$\lambda_F \leq \sum_{\Phi(F') < \Phi(F)}  \left( \frac{1}{2} \dim(F') + 1 \right) ,$$
where the sum is over all fixed components $F'$ such that $\Phi(F') < \Phi(F)$.
\end{Lemma}

\begin{proof}
Pick a fixed component $F$.
Let $N = \sum_{\Phi(F') < \Phi(F)} \left(\frac{1}{2} \dim(F')+ 1\right)$. 
By Lemma~\ref{le:extend}, there exists an equivariant extension $w \in H^2_{S^1}(M;\R)$  of $\omega$
so that $w|_{F'} = [\omega|_{F'}] - \Phi(F') t$ for all fixed components $F'$.
Define $\beta \in H^{2 N}_{S^1}(M;\R)$ by
$$\beta = \prod_{\Phi(F') < \Phi(F)} (w + \Phi(F') t)^{\frac{1}{2}\dim(F')+ 1},$$
where now the product is over all fixed components $F'$ such that $\Phi(F') < \Phi(F)$.
Given any fixed component $F'$  
$$ \left((w + \Phi(F') t)|_{F'}\right)^{\frac{1}{2} \dim(F')+ 1}   = 
 \left[\omega|_{F'}\right]^{\frac{1}{2} \dim(F') + 1} = 0.$$
Hence, the restriction $\beta|_{F'}$ vanishes for all fixed components $F'$ such
that $\Phi(F') < \Phi(F)$.
In contrast, as a polynomial in $t$,
\begin{align*}
\beta|_F 
&= \prod_{\Phi(F') < \Phi(F)} \left( (w + \Phi(F') t)|_F \right)^{\frac{1}{2}\dim(F')+ 1}    \\
&= \prod_{\Phi(F') < \Phi(F)} \left([\omega|_F]   + \Phi(F)t - \Phi(F')t\right)
^{\frac{1}{2}\dim(F')+1}   \\
&= \prod_{\Phi(F') < \Phi(F)} (\Phi(F) t - \Phi(F') t)^{\frac{1}{2}\dim(F') + 1}    + \mbox{lower order terms}. 
\end{align*}
Hence $\beta|_F \neq 0$.  
Since $\beta|_F$ is a multiple of $e_{S^1}(N^-_F)$
by Corollary~\ref{cor:basis}, this implies that $N \geq \lambda_F$.
\end{proof}

\begin{Remark}\label{rk:order}
In fact, by eliminating some of the redundant factors of $\beta$, one can
show that there exist fixed components $F_1,\dots,F_k $
such that $\Phi(F_i) < \Phi(F)$ for all $i$, 
$\lambda_F \leq \sum_i  \left( \frac{1}{2} \dim(F_i) + 1 \right),$
{\em and}
$\Phi(F_i) \neq \Phi(F_j)$ for all $i \neq j$.
\end{Remark}

When  $H^{2i}(M;\R) = H^{2i}(\CP^n;\R)$ for all $i$, there are
not many fixed components, and so this determines the order of
the fixed components under the moment map.

\begin{Lemma}\label{oneeach}
Let the circle act  on a
compact symplectic  manifold $(M,\omega)$
with moment map $\Phi \colon M \to \R$; assume that
$H^{2i}(M;\R)   = H^{2i}(\CP^n;\R)$ for all $i$.
There exists a unique fixed component $F$
such that $2 \lambda_F \leq 2i \leq 2 \lambda_F + \dim(F)$
for all $i$ such that $0 \leq 2i \leq 2n$.
\end{Lemma}

\begin{proof}
Since every fixed component $F$ is symplectic,
$H^{2i}(F;\R) \neq 0$
for every integer $i$ such that $0 \leq 2i \leq \dim (F)$.
The first claim  thus follows immediately from the assumption
and the fact that
moment maps are perfect Morse-Bott functions.
\end{proof}

\begin{Proposition}\label{cor:order}
Let the circle act on a compact symplectic manifold $(M,\omega)$
with moment map $\Phi \colon M \to \R$; assume that $H^{2i}(M;\R) =
H^{2i}(\CP^n;\R)$ for all $i$.
Then 
for all fixed components $F$ and $F'$,
$$
\Phi(F') < \Phi(F)
\quad \mbox{exactly if} \quad 
\lambda_{F'} < \lambda_F 
.$$
\end{Proposition}

\begin{proof}
Consider any fixed component $F$.
Define
\begin{equation}\label{o1}
N = \sum_{\Phi(F') < \Phi(F)} \left (\frac{1}{2}\dim(F') + 1 \right),
\end{equation}
where the sum is over all fixed components $F'$ such that $\Phi(F')< \Phi(F)$.
Since $H^{2i}(M;\R) = H^{2i}(\CP^n;\R)$ for all $i$,
Lemma~\ref{oneeach}  implies that 
there exists a unique fixed component $F'$ such that $2 \lambda_{F'} \leq 2i \leq 2
\lambda_{F'} + \dim(F')$
for all $i$ such that $0 \leq 2i \leq 2n$.
Therefore, since the fixed components are all even dimensional,
\begin{equation}\label{o3}
\sum_{\lambda_{F'} < N} \left( \frac{1}{2}\dim(F')+ 1\right)  \leq N,
\end{equation}
where here the sum is over all fixed components $F'$ such that $\lambda_{F'}  < N$.
On the other hand, 
by Lemma~\ref{le:order}, for any fixed component $F'$ with
$\Phi(F') \leq \Phi(F)$,
\begin{equation*}
\lambda_{F'} \leq \sum_{\Phi(F'') <  \Phi(F')} \left (\frac{1}{2}\dim(F'') + 1 \right)
 \leq \sum_{\Phi(F'') <  \Phi(F)} \left (\frac{1}{2}\dim(F'') + 1 \right)
= N ,
\end{equation*}
with equality impossible unless $\Phi(F')=\Phi(F)$.
In particular, $\lambda_{F'}< N$ for all fixed components
$F'$ such that $\Phi(F')<\Phi(F)$, and so equations \eqref{o1} and \eqref{o3} imply
that $\lambda_{F'} <N$ exactly if $\Phi(F') < \Phi(F)$.
Since $\lambda_F\leq N$, we can conclude that $\lambda_F=N$;
this proves the claim.
\end{proof}

For our main proposition, we will also need the fact
that the moment map is also a perfect Morse-Bott function over any field
whenever $H^{2i}(M;\R) = H^{2i}(\CP^n;\R)$ for all $i$,

\begin{Lemma}\label{Zperfect}
Let the circle act  on a
compact symplectic  manifold $(M,\omega)$
with moment map $\Phi \colon M \to \R$; assume that
$H^{2i}(M;\R)   = H^{2i}(\CP^n;\R)$ for all $i$.
Then
$$ H^j(M;\Z) =  \bigoplus_{F \subset M^{S^1}}  H^{j - 2\lambda_F}(F;\Z) \quad
\forall \ j,$$
where the sum is over all fixed components.
\end{Lemma}

\begin{proof}
Consider a fixed component $F$ and integer $j$
so that $H^{j- 2 \lambda_F}(F;\Z) \neq 0$.
Let $M^\pm = \Phi^{-1}(-\infty, \Phi(F) \pm \epsilon)$ for 
$\epsilon > 0$ sufficiently small.
Since each fixed component is even dimensional,  Lemma~\ref{oneeach} 
implies that
$H^{j- 2 \lambda_{F'}}(F';\Z) = H^{j-2 \lambda_{F'} + 1}(F';\Z)= 0$
for all fixed components $F' \neq F$.
Therefore, $H^{j}(M^-;\Z) = H^{j+ 1}(M^-;\Z) = 0,$
and so $H^j(M^+;\Z) = H^{j - 2\lambda_F}(F;\Z)$.
By a similar argument, $H^j(M^+;\Z)  = H^j(M;\Z)$.

\end{proof}

Given any fixed component $F$, let  $\Gamma_F \in \Z$
denote the sum of the weights at $F$.

\begin{Lemma}\labell{c1}
Let the circle  act  on a
compact symplectic  manifold $(M,\omega)$
with moment map $\Phi \colon M \to \R$.
Then $c_1(M) \neq 0$.
\end{Lemma}

\begin{proof}
Given a point $p$ in a fixed component $F$,
$c^{S^1}_1(M)|_p = \Gamma_F \, t$.
Each weight at the minimal fixed component  $F$ is positive,
while each weight at the maximal fixed component $F'$ is negative.
Hence, $\Gamma_{F} \neq \Gamma_{F'}$;
this implies that $c^{S^1}_1(M)$ is not a multiple of $t$.
\end{proof}

Finally, we state our main proposition.

\begin{Proposition}\labell{prop:isom}
Let the circle act  on a compact symplectic manifold $(M,\omega)$
with moment map $\Phi \colon M \to \R$;
assume that $H^j(M;\Z) = H^j(\CP^n;\Z)$ for all $j$.

For each integer $i$ such that $0 \leq 2i \leq 2n$,
there exists a unique fixed component $F_i$ 
so that $H^{2i - 2 \lambda_{F_i}}(F_i;\Z) = \Z$; 
let $u_i$ be a generator of $H^{2i - 2 \lambda_{F_i}}(F_i;\Z)$, and
let $c_1(M)^{i - \lambda_{F_i}}|_{F_i}  =  m_i  u_i$.
Then the cohomology class  $$\alpha_i =
 \frac{\Lambda_{F_i}^- }{m_i}
 \left(c^{S^1}_1(M) -  \Gamma_{F_i}\, t \right)^{i - \lambda_{F_i}} 
\prod_{\lambda_{F'} < \lambda_{F_i}} 
\left( \frac{c^{S^1}_1(M) -  \Gamma_{F'}\, t } {\Gamma_{F_i} - \Gamma_{F'} } \right)
^{\frac{1}{2}\dim (F') + 1}
$$
is well-defined and lies in $H_{S^1}^{2i}(M;\Z)$.
(Here, the product is over all fixed components $F'$ such that $\lambda_{F'} < 
\lambda_{F_i}$.)
Moreover, the classes $\alpha_0,\ldots,\alpha_n$ form a basis for
$H^*_{S^1}(M;\Z)$ as a $H^*(\CP^\infty;\Z) = \Z[t]$ module.
\end{Proposition}

\begin{proof}
By Lemma~\ref{oneeach}, for each integer $i$ such that
$0 \leq 2i \leq 2n$ there exists a unique fixed component
$F_i$ such that $2 \lambda_{F_i} \leq 2i \leq 2\lambda_{F_i}+\dim(F_i)$.
By Lemma~\ref{Zperfect}, 
$H^{2i - 2 \lambda_{F_i}}(F_i;\Z) \simeq \Z$; let $u_i$ be a generator.
Then  $u_0,\dots,u_n$ is a basis for $H^*(M^{S^1})$.
By  Proposition~\ref{basis}, 
for each such $i$ there exists a class $\alpha_i \in H^{2i}_{S^1}(M;\Z)$ so that
\begin{enumerate}
\item $\alpha_i|_{F_i} =  u_i \,  e_{S^1}(N^-_{F_i})$, and
\item $\alpha_i|_{F'} = 0$ for all other fixed components 
$F'$ with $\Phi(F' ) \leq \Phi(F_i)$.
\end{enumerate}
Moreover, 
$\alpha_0,\dots,\alpha_n$ is a basis for $H^*_{S^1}(M;\Z)$ as
an $H^*(\CP^\infty;\Z)$ module.

Since each fixed component has even dimension,
Lemma~\ref{oneeach} implies that  for each fixed component $F$,
$$
 \sum_{\lambda_{F'} < \lambda_{F}} \left( \frac{1}{2} \dim(F') + 1 \right) =
 \lambda_{F} ,$$ 
where here the sum is over all fixed components $F'$ such that $\lambda_{F'} < \lambda_F$.
Hence, 
for each $i$ such that $0\leq 2i \leq \dim(M)$
we may define
$$\beta_i =
\left(c^{S^1}_1(M) -  \Gamma_{F_i} \, t \right)^{i - \lambda_{F_i}}
\prod_{\lambda_{F'} < \lambda_{F_i}} 
\bigl( {c^{S^1}_1(M) -  \Gamma_{F'} \, t } \bigr)
^{\frac{1}{2}\dim(F') + 1} \ \in H^{2i}_{S^1}(M;\Z).$$

Moreover, for every fixed component $F'$, 
$$(c_1^{S^1}(M) - \Gamma_{F'} t)^{\frac{1}{2} \dim(F')  +1}|_{F'}
= \left(c_1(M)|_{F'}\right)^{\frac{1}{2} \dim (F') + 1} = 0.$$
Therefore, the  restriction $\beta_i|_{F'}$ vanishes for every
fixed component $F'$ such that $\lambda_{F'}  < \lambda_{F_i}$.
Combining  Corollary~\ref{cor:basis} and Proposition~\ref{cor:order}, 
this implies that 
$$\beta_i = \sum_{\lambda_{F_j} \geq \lambda_{F_i}} x_j \alpha_j,$$
where $x_j \in H^*(\CP^\infty;\Z) = \Z[t]$ for all $j$.
Here, the sum is over $j$ such that $\lambda_{F_j} \geq \lambda_{F_i}$.
Since $\beta_i \in H^{2i}_{S^1}(M;\Z)$, Lemma~\ref{oneeach} and degree considerations
imply that
\begin{equation}\labell{isom1}
\beta_i = \sum_{\substack{F_j = F_i \\ j \leq i}} x_j \alpha_j,
\end{equation}
where now the sum is over $j$ such that $F_j = F_i$ and $j \leq i$.

Additionally,
recall that  $c_1(M)^{i - \lambda_{F_i}}|_{F_i} = m_i u_i$.
As a polynomial in $t$,
\begin{equation}\labell{isom2}
\begin{split}
\beta_i|_{F_i}  &= 
c_1(M)^{i - \lambda_{F_i}}  |_{F_i}
  \prod_{\lambda_{F'} < \lambda_{F}} 
\bigl(
c_1(M)|_F +  (\Gamma_F  -  \Gamma_{F'} ) t 
\bigr)
^{\frac{1}{2}\dim(F') + 1} \\
& = m_i  u_i \,  \prod_{\lambda_{F'} < \lambda_{F_i}}
(\Gamma_{F_i} - \Gamma_{F'})^{\frac{1}{2}\dim(F') + 1} \, t^{\lambda_{F_i}}
+ \mbox{lower order terms}.
\end{split}
\end{equation}
On the other hand,   for any $j$ such that $F_j = F_i$,
\begin{equation}\labell{isom3}
\alpha_j|_{F_i} =  u_j e_{S^1}(N_{F_i}^-) = u_j \Lambda_{F_i}^-t^{\lambda_{F_i}}  + \mbox{lower order terms}. 
\end{equation}
Comparing  equations  \eqref{isom1}, \eqref{isom2}, and \eqref{isom3},
we see that 
$$
\beta_i = 
\frac{m_i}{\Lambda^-_{F_i}}
\prod_{\lambda_{F'} < \lambda_{F_i}} 
{(\Gamma_{F_i} - \Gamma_{F'})^{ \frac{1}{2}\dim(F') + 1} } \alpha_i.
$$

Finally, by Lemma~\ref{c1}, $c_1(M) \neq 0$.
Since $H^2(M;\R) = \R$, this implies that $c_1(M)$ is
a non-zero multiple of $\omega$.   Therefore, $c_1(M)^i \neq 0$.
Since $\beta_i$ maps to $c_1(M)^i$ under the natural
restriction map from the equivariant cohomology of $M$
to the ordinary cohomology, this clearly implies that $\beta_i \neq 0$.
\end{proof}

We can use this proposition (and Remark~\ref{ABBV}) to obtain
particularly nice description of  
the equivariant cohomology and total Chern class
in the following special cases. 

%% used in intro

\begin{Corollary}\label{cor:eqisom}
Let the circle act  on compact symplectic manifolds $(M,\omega)$ and 
$(\Hat{M},\Hat{\omega})$ with moment maps $\Phi \colon M \to \R$ and $\Hat\Phi \colon
\Hat{M} \to \R$, respectively;
assume that $H^j(\Hat{M};\Z) = H^j(\CP^n;\Z)$ for all $j$.
If $f \colon M^{S^1} \to \Hat{M}^{S^1}$
is a diffeomorphism and
$T(M)\big|_{M^{S^1}} \cong 
f^*\bigl( T (\Hat{M})  \big|_{  \Hat{M} ^{S^1} } \bigr), $ 
then
$H_{S^1}^*(M;\Z)\big|_{M^{S^1}} =
f^*\bigl(H_{S^1}^*(\Hat{M} ;\Z)\big|_{\Hat{M}^{S^1}} \bigr)
$ and 
$
c(M)\big|_{M^{S^1}}
= 
f^*\bigl(c(\Hat{M})\big|_{\Hat{M}^{S^1}} \bigr) 
.$
\end{Corollary}

\begin{proof}
Lemma~\ref{Zperfect} implies that
$$\bigoplus_{\Hat{F} \subset \Hat{M}^{S^1}} H^{j - 2\lambda_{\Hat{F}}}(\Hat{F};\Z)
 = H^j(\Hat{M};\Z) = H^j(\CP^n;\Z) \quad \forall \ j.$$
Since, $T(M)\big|_{M^{S^1}} \cong 
f^*\bigl( T (\Hat{M})  \big|_{  \Hat{M} ^{S^1} } \bigr), $ 
this implies that 
$$\bigoplus_{F \subset M^{S^1}} H^{j - 2\lambda_F}(F;\Z)
= H^j(\CP^n;\Z) \quad \forall \ j. $$
As in Lemma~\ref{Zperfect}, this implies that
$$H^j(M;\Z) = H^j(\CP^n;\Z) \quad \forall\ j.$$
The claim now follows immediately from Proposition~\ref{prop:isom}.
\end{proof}

%%% used in intro

\begin{Corollary}\label{cor:eqisol}
Let the circle act  on compact symplectic manifold $(M,\omega)$ with
moment map $\Phi \colon M \to \ft^*$;  
assume that there is a unique fixed point of index $2i$ for
all $i$ such that $0 \leq i \leq n$.
As a $H^*(\CP^\infty;\Z) = \Z[t]$ module,
$H_{S^1}^*(M;\Z)$ 
is freely generated by
$1,\alpha_1,\ldots,\alpha_n$, where 
$$\alpha_i =    \Lambda^-_{p_i} \, \prod_{j = 0}^{i - 1}
\frac {c^{S^1}_1(M)  - \Gamma_{p_j} t} {\Gamma_{p_i}  - \Gamma_{p_j}}.$$
(In particular,  $\alpha_i \in H_{S^1}^{2i}(M;\Z)$ for all $i$.)
\end{Corollary}

%% used in intro

\begin{Example}\label{ex:eq}
Let the circle act on a compact symplectic manifold $(M,\omega)$ with
moment map $\Phi \colon M \to \ft^*$.
Assume that the fixed set consists of four points
$p_0,p_1,p_2$ and $p_3$
with weights $\{1,2,3\}, \{1,-1,l\}, \{1,-1,-l\},$
and $\{-1,-2,-3\}$, respectively.  
As a $H^*(\CP^\infty;\Z) = \Z[t]$ module, $H_{S^1}^*(M;\Z)$ is freely
generated by $1, \alpha_1,\alpha_2,\alpha_3$, where
\begin{gather*}
\alpha_1|_{p_1} = t, \
\alpha_1|_{p_2} = \frac{6+ l}{6 - l}t , \
\alpha_1|_{p_3} = \frac{12}{6 - l} t, \
\alpha_2|_{p_2} =   l t^2,  \\
\alpha_2|_{p_3} = 6 t^2,  \
 \alpha_3|_{p_3} = 6 t^3,
\ \mbox{and} \ \alpha_i|_{p_j} = 0 \ \forall \ j < i;
\ \mbox{moreover}, \\
c^{S^1}(M)|_{p_0} =  1 + 6t + 11 t^2 + 6 t^3, \ \ 
c^{S^1}(M)|_{p_1} =  1 + l t  -  t^2 - l t^3,  \\
c^{S^1}(M)|_{p_2} =  1 - lt  -  t^2 + l t^3, \ \  \mbox{and} \ \ 
c^{S^1}(M)|_{p_3} =  1 - 6t + 11 t^2 - 6 t^3. \\
\end{gather*}
\end{Example}

In the above cases, the fixed set is torsion free and so
formality holds.
Thus, the ordinary cohomology $H^*(M;\Z)$ and total Chern class $c(M)$
also very easy to describe.

%% used in intro

\begin{Corollary}\label{cor:isom}
Let the circle act  on compact symplectic manifolds $(M,\omega)$ and 
$(\Hat{M},\Hat{\omega})$ with moment maps $\Phi \colon M \to \R$
and $\Hat{\Phi} \colon \Hat{M} \to \R$, respectively;
assume that $H^j(\Hat{M};\Z) = H^j(\CP^n;\Z)$ for all $j$.
If $f \colon M^{S^1} \to \Hat{M}^{S^1}$ is a diffeomorphism  and
$T(M)\big|_{M^{S^1}} \cong 
f^*\bigl( T (\Hat{M})  \big|_{  \Hat{M} ^{S^1} } \bigr), $ 
then 
$f$ induces an isomorphism
$f^\sharp \colon H^*(\Hat{M};\Z) \to H^*(M;\Z)$ so that 
$
c(M)
= 
f^\sharp\bigl(c(\Hat{M})\bigr). 
$
\end{Corollary}

\begin{Remark}\labell{rmk:cob}
More generally,
let the circle act  on compact symplectic manifolds $(M,\omega)$ and 
$(\Hat{M},\Hat{\omega})$ with moment maps $\Phi \colon M \to \R$ and $\Hat\Phi \colon
\Hat{M} \to \R$, respectively;
let $f \colon M^{S^1} \to \Hat{M}^{S^1}$
be an orientation preserving  diffeomorphism such that
$T(M)\big|_{M^{S^1}} \cong 
f^*\bigl( T (\Hat{M})  \big|_{  \Hat{M} ^{S^1} } \bigr).$ 

Since moment maps are (equivariantly) perfect Morse-Bott functions,
$M$ and $\Hat{M}$ have the same  (equivariant) Betti numbers.
Moreover, Theorem~\ref{ABBV}
immediately implies that that $M$ and $\Hat{M}$ have the
same (equivariant) Chern {\em numbers}. 
Alternately, by \cite{GGK}, $M$ and $\Hat{M}$ are cobordant
as oriented equivariant stable-complex manifolds,
and so have the same (equivariant)
Chern numbers.  
Additionally, if $H_2(M;\R) = H_2(\Hat{M};\R)$,
then since  $c_1(M) \neq 0$ by Lemma~\ref{c1}, $c_1^{S^1}(M)$  and $t$ are a basis
for  $H^2_{S^1}(M;\R)$ as a vector space.
By assumption,
$$f^*(c^{S^1}(\Hat{M})\big|_{\Hat{M}^{S^1}}) = c^{S^1}(M)\big|_{M^{S^1}}.$$
Hence, after  possibly multiplying $\Hat{\omega}$ by a constant and adding a constant to
$\Hat{\Phi}$, $$f^*([\Hat{\omega} - \Hat{\Phi} t]\big|_{\Hat{M}^{S^1}} ) = [ \omega - \Phi t]\big|_{M^{S^1}}.$$
Therefore, by \cite{GGK}, $M$ and $\Hat{M}$ are also cobordant as Hamiltonian
$S^1$-spaces. 
This  (or  Theorem~\ref{ABBV}) implies that $M$ and $\Hat{M}$ have the same
Duistermaat-Heckman measure. More generally,  any product of equivariant Chern
classes and powers of the class $[\omega - \Phi t]$ have the same integral over $M$ and $\Hat{M}$.

However, in general cobordant manifolds do {\em not} have
isomorphic (equivariant) cohomology rings. 
For example, let $M$ be the blow-up of $\CP^2$ at $[0,1,0]$,
where $S^1$ acts on $\CP^2$ by $$\lambda \cdot [x_0,x_1,x_2] = [ x_0, \lambda x_1, \lambda^2 x_2],$$
and let $S^1$ act on $\Hat{M} = \CP^1 \times \CP^1$ by 
$$\lambda \cdot ([y_0,y_1], [z_0,z_1] ) = ([y_0, \lambda y_1], [z_0,\lambda^2 z_1]).$$
For the appropriate choice of symplectic forms,
$M$ and $\Hat{M}$  are cobordant as stable-complex Hamiltonian $S^1$-spaces,
but as rings $H^*(M;\Z) \not\cong H^*(\Hat{M};\Z).$

\end{Remark}

\begin{Remark}
Alternately, let the circle act  on compact symplectic manifolds $(M,\omega)$ and 
$(\Hat{M},\Hat{\omega})$ with moment maps $\Phi \colon M \to \R$
and $\Hat{\Phi} \colon \Hat{M} \to \R$, respectively;
assume that $H^j(\Hat{M};\Z) = H^j(\CP^n;\Z)$ for all $j$.
%Alternately, even if there is no diffeomorphism $f$, 
Also assume that there is a bijection from the fixed
components $F_1,\dots,F_k$ of $M$ to the fixed components $\Hat{F}_1,\dots,\Hat{F}_k$ of $\Hat{M}$
and that there exists an 
isomorphism $f^* \colon H^*_{S^1}(\Hat{F}_i;\Z) \to H^*_{S^1}(F_i;\Z)$
such that $f^*(c^{S^1}(\Hat{M})|_{\Hat{F}_i}) = c^{S^1}(M)|_{F_i}$ for all $i$.
Then all the symmetric polynomials in the weights at $F$ and $\Hat{F}$ are
the same, and so they have the same weights, and hence the same index.
Therefore, the conclusions of Corollaries~\ref{cor:eqisom} and 
~\ref{cor:eqisom} still hold.
\end{Remark}
%%S used in intro

\begin{Corollary}\labell{cor:isol}
Let the circle  act on a compact symplectic 
manifold $(M,\omega)$ with moment map $\Phi \colon M \to \ft^*$;
assume that there is a unique fixed point $p_i$
of index $2i$  for all $i$ such that $0 \leq 2i \leq 2n$.
As a group,  $H^*(M;\Z)$ is freely generated by
$1,\Tilde{\alpha}_1,  \ldots, \Tilde{\alpha}_n$,
where   
\begin{gather*}
\Tilde{\alpha}_i
= 
\frac
{\Lambda_{p_i}^- }
{(\Lambda_{p_1}^-)^i}
\, \frac
{ (\Gamma_{p_1} - \Gamma_{p_0})^i}
{ \prod_{j = 0}^{i-1} (\Gamma_{p_i} - \Gamma_{p_j})}
\;
(\Tilde{\alpha}_1)^i 
; 
\quad \mbox{moreover} \\
\begin{split}
c_i(M) & =
\frac{ \Lambda_{p_{i}}^+}{\prod_{j = i+1}^{n} 
(\Gamma_{p_{i}} - \Gamma_{p_j}) }
\left(  \sum_{k = 0}^i 
{c^{S^1}_i(M)|_{p_k}}
\frac
{\prod_{j = i+1}^n (\Gamma_{p_k} - \Gamma_{p_j})}
{ t^i \Lambda_{p_k}  } 
\right) \Tilde{\alpha}_i \\
& = \frac
{1}{ \Lambda_{p_{i}}^- t^i}
\prod_{j = 0}^{i-1} (\Gamma_{p_{i}} - \Gamma_{p_j}) 
\left(  \sum_{k = 0}^i 
\frac {c^{S^1}_i(M)|_{p_k}}
%\prod_{j \in \{0,\dots,i\} \setminus \{k\}} (\Gamma_{p_k} - \Gamma_{p_j})
{\prod_{j \in \{0,\dots,\hat{k},\dots i\} } (\Gamma_{p_k} - \Gamma_{p_j})}
\right) \Tilde{\alpha}_i. 
\end{split}
\end{gather*}
(In particular, $\Tilde{\alpha}_i \in H^{2i}(M;\Z)$ for  all $i$.)
\end{Corollary}

\begin{proof}
As a $H^*(\CP^\infty;\Z) = \Z[t]$ module,
$H_{S^1}^*(M;\Z)$ 
is freely generated by the classes $1,\alpha_1,\ldots,\alpha_n$, where 
$$\alpha_i =    \Lambda^-_{p_i} \, \prod_{j = 0}^{i - 1}
\frac {c^{S^1}_1(M)  - \Gamma_{p_j} t} {\Gamma_{p_i}  - \Gamma_{p_j}} 
\in H_{S^1}^{2i}(M;\Z) \quad  \forall \ i.$$
The image of $\alpha_i$ under the restriction map from equivariant cohomology to ordinary cohomology is
$$\Tilde{\alpha}_i  =  \Lambda_{p_i}^- \frac {c_1(M)^i}
{\prod_{j = 0}^{i-1} (\Gamma_{p_i} - \Gamma_{p_j}) } \in H^{2i}(M;\Z).
$$
By formality, 
$H^*(M;\Z)$ is freely generated (as a group) by
$\Tilde{\alpha}_0, \ldots, \Tilde{\alpha}_n$.

Applying Proposition~\ref{prop:isom} to the reversed circle action with moment map
$-\Phi$, we may define
$$\beta_{n-i} =   \Lambda^+_{p_{i}} \, \prod_{j = i+1}^{n}
\frac {c^{S^1}_1(M)  - \Gamma_{p_j} t} {\Gamma_{p_{i}}  - \Gamma_{p_j}} 
\in H_{S^1}^{2n - 2i}(M;\Z);$$
let $\Tilde{\beta}_{n-i} \in H^{2n-2i}(M;\Z)$ be the image of $\beta_{n-i}$ under the restriction
map from equivariant cohomology to ordinary cohomology.
Note that $$
( \alpha_i \beta_{n-i} )|_{p_k} = 
\begin{cases}
\Lambda_{p_i} & k = i  \\
0 & \mbox{otherwise}. \\
\end{cases} $$
Therefore, Theorem~\ref{ABBV} implies that
\begin{gather}
\int_M \alpha_i \beta_{n-i} = 1, \quad \mbox{and so} \label{dodo}\\
\int_M \Tilde{\alpha}_i \Tilde{\beta}_{n-i} = 1.
\end{gather}
Since  $c_i(M)$  is a multiple of $\Tilde{\alpha}_i$, this
implies that
$$c_i(M) = \left( \int_M c_i(M) \Tilde{\beta}_{n-i} \right) \Tilde{\alpha}_i = \left( \int_M  c_i^{S^1}(M) \beta_{n-i} \right)
\Tilde{\alpha}_i.$$
Since  
\begin{equation*}
\beta_{n-i}|_{p_k} =
\begin{cases}
 \Lambda_{p_i}^+ t^{n-i} \prod_{j=i+1}^n \frac
{\Gamma_{p_k} - \Gamma_{p_j}}
{\Gamma_{p_{i}} - \Gamma_{p_j}}
& \mbox{if}\  0  \leq k \leq i, \\
0  & \mbox{if} \  i <   k  \leq  n,  
\end{cases}
\end{equation*}
the first equality follows immediately  from  Theorem~\ref{ABBV}.

Finally, since $\int_M \alpha_k \beta_{n-k} = 1$ for all $k \in \{0,\dots,n\}$ by \eqref{dodo},
$$ \frac{\prod_{j \neq k}(\Gamma_{p_k} - \Gamma_{p_j}) }{\Lambda_{p_k}} =
\int_M \prod_{j \neq k} \left( c_1^{S^1}(M) - \Gamma_{p_j} \right)
=  \int_M c_1(M)^n.
$$
Hence, 
$$ \frac{\prod_{j \neq k}(\Gamma_{p_k} - \Gamma_{p_j}) }{\Lambda_{p_k}} 
= \frac{\prod_{j \neq l}(\Gamma_{p_l} - \Gamma_{p_j}) }{\Lambda_{p_l}} \quad \forall
k,l \in \{0,\dots,n\}.$$
The second equality follows immediately.
\end{proof}

Note that the equations for $c_i(M)$ above simplify in some cases. 
Not only is $c_0(M) = 1$ and $c_n(M) = (n +1) \Tilde{\alpha}_n$, but also
\begin{gather*}
%\begin{equation}
%\begin{split}
 c_1(M) = \frac{\Gamma_{p_1} - \Gamma_{p_0}}{\Lambda^-_{p_1}} \Tilde{\alpha}_1, \quad \mbox{and}  \\
 c_2(M) = \frac
{{c^{S^1}_2\! (M)|_{p_0}}(\Gamma_{p_2} - \Gamma_{p_1})
- {c^{S^1}_2 \! (M)|_{p_1}}(\Gamma_{p_2} - \Gamma_{p_0})
+ {c^{S^1}_2 \!(M)|_{p_2}}(\Gamma_{p_1} - \Gamma_{p_0})}
{(\Gamma_{p_1} - \Gamma_{p_0}) \Lambda^-_{p_2} t^2} 
\Tilde{\alpha}_2.
%\end{split}
%\end{equation}
\end{gather*}

\begin{Example}\labell{ex:non}
Let the circle act on a symplectic manifold $(M,\omega)$ with
moment map $\Phi \colon M \to \ft^*$.
Assume that the fixed set consists of four points
$p_0,p_1,p_2$ and $p_3$
with weights $\{1,2,3\}, \{1,-1,l\}, \{1,-1,-l\},$
and $\{-1,-2,-3\}$, respectively.  
Then
$
\Lambda^-_{p_1} = -1,\
\Lambda^-_{p_2} = l,  \
\Gamma_{p_0} =  6,\
\Gamma_{p_1} = -\Gamma_{p_2} = l, \
c_2(M)|_{p_0} = 11t^2,\
$ and $
c_2(M)|_{p_1} = c_2(M)|_{p_2} = - t^2$.
Therefore, Corollary~\ref{cor:isol} (and the equations above)
immediately imply that
\begin{align*}
H^*(M;\Z) &= \Z[x,y]/\Bigl(x^2 - \frac{2(l + 6)}{(6 - l)^2} y,y^2 \Bigr), \quad \mbox{and}\\
c(M) & = 1 + (6 - l) x + \frac{24}{6 - l} y + 4 xy,
\end{align*}
where $x$ has degree $2$ and $y$ has degree $4$.
\end{Example}

The final lemma gives the relationship between  the moment
image and the sum of the weights at fixed components.

\begin{Lemma}\labell{comp}
Let the circle act  on a
compact symplectic manifold $(M,\omega)$
with moment map $\Phi \colon M \to \R$.
If $H^{2}(M;\R)   = \R$, then
for all fixed components $F$ and $F'$,
$$
\Gamma_{F} > \Gamma_{F'}
\quad \mbox{exactly if} \quad 
\Phi(F) < \Phi(F') 
.$$
\end{Lemma}

\begin{proof}
By Lemma~\ref{le:extend}, there exists an equivariant extension
$w \in H^2_{S^1}(M;\R)$ of $[\omega]$ so that
$w|_p = - \Phi(p) t$ for all fixed points $p$.
Since $H^2(M;\R) = \R$ and $ [\omega] \neq 0$,
$c^{S^1}_1(M) = a w   + b t$  for some real numbers $a$ and $b$. 
Therefore, given any fixed component $F$,
$\Gamma_F \, t = c^{S^1}_1(M)|_p = (aw + b t)|_p 
= (-a \Phi(F) + b) t$
for all $p \in F$.
On the other hand, if $F$ is the minimal fixed component and
$F'$ is the maximal fixed component then 
$\Phi(F) < \Phi(F')$ and  $\Gamma_F > 0 > \Gamma_{F'}$.
Therefore, $a > 0$.
\end{proof}

\section{The case that the fixed set is not discrete}
\label{s:notdiscrete}

We now return to the $6$-dimensional case.
In this section, we prove Theorem~\ref{thm:2} in the case that the fixed
set is not discrete. In fact, in this case only the first two
possibilities can arise.

\begin{Proposition}\label{prop:notdiscrete}
Let the circle act faithfully on a $6$-dimensional 
compact symplectic manifold 
$(M,\omega)$ with moment map $\Phi \colon M \to \R$.
Assume that $H^2(M,\R) = \R$  
and the fixed set is not discrete.
Then one of the following two statements is true:
\begin{itemize}
\item [(A)]
There is a subgroup $S^1 \subset \SU(4)$
and an orientation preserving diffeomorphism
$f \colon M^{S^1} \to \bigl(\CP^3 \bigr)^{S^1}$
so that
$T (M)|_{M^{S^1}} \cong f^*\left(
T \bigl(\CP^3 \bigr) \big|_{ ( \CP^3 )^{S^1} }
\right). $ 
\item [(B)]
There is a subgroup $S^1 \subset \SO(5)$
and an orientation preserving diffeomorphism
$f \colon M^{S^1} \to \Tilde{G}_2(\R^5)^{S^1} $  
so that 
$T(M)|_{M^{S^1}} \cong 
f^*\left(T \bigl(\Tilde{G}_2(\R^5) \bigr) 
\big|_{\Tilde{G}_2(\R^5)^{S^1}}\right) $. 
\end{itemize}
\end{Proposition}

\begin{Remark}
The manifolds described above do not contain an isolated fixed
point with three distinct weights.  
A fortiori, they do not contain a pair of isotropy spheres
which intersect in two points.
\end{Remark}

\begin{Remark}\labell{rmk:cob2}
In fact, after possibly multiplying 
the standard symplectic form $\Hat{\omega}$ on $\CP^3$ (or $\Tilde{G}_2(\R^5)$) 
by a constant, we may assume that $f$ is a symplectomorphism on each fixed component.
To see this, first note that the argument in Remark~\ref{rmk:cob} 
implies that $f^*([\Hat{\omega}]|_{\Hat{M}^{S^1}}) = [\omega]|_{M^{S^1}}$.
If each fixed component is at most $2$-dimensional,  this immediately implies
that there exists a symplectomorphism $f' \colon M^{S^1} \to \Hat{M}^{S^1}$
which is homotopic to $f$.  Otherwise, the claim follows
from \cite{Del}; see Case IV below.
\end{Remark}

We will need the following lemma to analyze the isotropy
submanifolds which might arise.

\begin{Lemma}\labell{fourdimb}
Let the circle act on a 
$4$-dimensional compact symplectic manifold $(Z,\sigma)$
with moment map $\Psi \colon Z \to \R$.
\begin{enumerate}
\item If $Z^{S^1}$ consists of a minimal surface $\Sigma$ and one point, 
then $e(N_\Sigma)$, the  Euler class of the normal bundle to $\Sigma$, is 
the positive generator of $H^2(\Sigma;\Z)$.

\item If $Z^{S^1}$ consists of a minimal surface $\Sigma$ and two points, 
then $e(N_\Sigma) = 0$.
\item
If $Z^{S^1}$ consists of a minimal surfaces $\Sigma$ and a maximal surface $\Sigma'$,
then there exists $a \in \Z$ so that
 $e(N_\Sigma) =  au$ and
$e(N_{\Sigma'}) = - au'$.
Here, $u$ and $u'$ are the positive generators of $H^2(\Sigma;\Z)$,
and $H^2(\Sigma';\Z)$, respectively. 
\end{enumerate}
\end{Lemma}

\begin{proof}
By dividing out by a finite subgroup, we may assume that the action
is faithful.  Let $e_{S^1}(N_\Sigma) = t + e(N_\Sigma) = t + au$, where
$H^*(\CP^\infty;\Z) = \Z[t]$.
Recall that every isotropy sphere contains
exactly two isolated fixed points.

In case (1), since $\Psi$ is a Morse-Bott function the fixed point must have index $4$.
Since there is only one isolated fixed point,
there are no isotropy spheres. Hence,  both weights at the isolated fixed  point are $-1$.
Since the degree of $1$ is $0 < 4$, $\int_Z 1 = 0$. Hence
Theorem~\ref{ABBV} and Remark~\ref{rmk:ABBV} imply   that
$ -a + 1 = 0$, that is, $a = 1$. 

In case (2), since $\Psi$ is a Morse-Bott function one fixed point must have index $2$
and one must have index $4$.
There are no isotropy spheres except possibly one joining
these two points.
Therefore, the  weights at these points
are  $\{-1,l\}$ and  $\{-1,-l\}$
for some natural number $l$. 
Since $\int_Z 1 = 0$, Theorem~\ref{ABBV} implies  that
$- a  - \frac{1}{l} + \frac{1}{l} = 0$, that is, $a = 0$.

In case (3),
let $e_{S^1}(N_{\Sigma'}) = -t + e(N_{\Sigma'}) = -t + bu'$.
Since $\int_Z 1 = 0$, Theorem~\ref{ABBV} implies  that
$- a  - b = 0$. 
\end{proof}

We will spend the remainder of this section proving
Proposition~\ref{prop:notdiscrete}. Let the circle
act on a symplectic manifold $(M,\omega)$ with moment map $\Phi \colon M \to \R$;
assume that  conditions of the proposition are satisfied.
By Remark~\ref{rmk:betti}, $H^{2i}(M;\R) = \R$ for all $i$ such that $0 \leq 2i \leq 6$.

Let  $N \subset  M^{\Z_k}$ be any $4$-dimensional 
isotropy submanifold.
Since  $\Phi|_N$ is a Morse-Bott function and the critical sets of $\Phi|_N$ are the fixed sets,
$N^{S^1}$ must contain at least two fixed components,
and at least one component must have index $0$ or $2$.

Finally, notice that the action
obtained by  reversing the circle action
and replacing $\Phi$ by $-\Phi$
also satisfies
the assumptions of Proposition~\ref{prop:notdiscrete}.
Moreover, if this new action satisfies the conclusions
of the proposition, then the original action does as well.
Therefore, we can replace $\Phi$ by $-\Phi$  at any time.
Given this symmetry, we only need to consider four cases.

\noindent
{\bf Case I: The fixed set is discrete except one minimal fixed surface.}

Since $\Phi$ is a perfect Morse-Bott function,
Poincar\'e duality implies that
the fixed set
consists of a minimal fixed sphere $\Sigma_0$ and two fixed
points $p_2$ and $p_3$ of index $4$ and $6$, respectively.
Since every vector bundle over a surface splits,
the normal bundle to $\Sigma_0$ is the direct sum of two
line bundles with equivariant Euler classes 
$mt + au$ and $n t + bu$, where $m$ and $n$ are natural numbers,
$a$ and $b$ are integers,
and $u$ is the positive generator of $H^2(\Sigma_0;\Z)$.

Assume first  that there is a $4$-dimensional 
isotropy submanifold  $N \subset M^{\Z_m}$
which contains both $p_2$ and $p_3$.
Since the intersection of any two
distinct $4$-dimensional isotropy submanifolds is  $2$-dimensional,
each component of the intersection must be either a fixed surface
or an isotropy sphere.  
Hence, every $4$-dimensional isotropy submanifold must contain
both fixed points.
So  there exists a natural number $l$ which is a multiple
of both  $m$ and $n$ so
that the weights at $p_2$ and $p_3$ are $\{-m,-n,l\}$
and  $\{-m,-n,-l\}$, respectively.
Applying  Lemma~\ref{fourdimb} to
the isotropy submanifold $N \subset M^{\Z_m}$, we see that $a = 0$.
Moreover, Lemma~\ref{modulo}
implies that  $2n = 0 \mod m$. Since $m \neq 1$ and $n$ and $m$ are relatively prime, 
this implies that $m = 2$ and $n = 1$.
Since the degree of $1$ is $0 < 6$ and the degree of $c^{S^1}_1(M)$ is  $2 < 6$,
$\int_M 1 = 0$ and $\int_M c^{S^1}_1(M) = 0$.
Therefore, Theorem~{\ref{ABBV} and Remark~\ref{rmk:ABBV} imply  that
$$ -\frac{b}{2} + \frac{1}{2l} - \frac{1}{2l} = 0 \quad \mbox{and} \quad
(1 - b) + \frac{l - 3}{2l} + \frac{l + 3}{2l} = 0.$$
These equations simplify to  $b = 0$ and $b = 2$,
which is impossible.

Therefore,  there is no $4$-dimensional isotropy submanifold
which contains both $p_2$ and $p_3$. 
If $m \neq 1$, then this implies that $M^{\Z_m}$ contains a 
$4$-dimensional
isotropy submanifold $N$ which contains $\Sigma_0$ and exactly
one fixed point.  By Lemma~\ref{fourdimb}, this implies
that $a = 1$.  A nearly identical argument
implies that if $n \neq 1$ then $b = 1$.
%Then every $4$-dimensional isotropy submanifold $N$ 
%must contain $\Sigma_0$ and  one fixed point. 
%By Lemma~\ref{fourdimb}, this implies that if
%$m \neq 1$ then $a = 1$, and  if $n \neq 1$ then $b = 1$.
Moreover, if there exist any isotropy spheres, then each 
one must contain  $p_2$ and $p_3$.
Hence, there exists a natural number $l$ so that the weights
at $p_2$ and $p_3$ are 
$\{-m, -m, l\}$  and 
$\{-n,-n,-l\}$, respectively.
Since $\int_M 1 = 0$ 
and $\int_M c^{S^1}_1(M)  = 0$, 
 Theorem~\ref{ABBV} and Remark~\ref{rmk:ABBV} imply  that
\begin{gather}
\labell{1Ia}
-\left(\frac{a}{m^2 n} + \frac{b}{m n^2}\right) 
+ \frac{1}{l m^2} - \frac{1}{l n^2} = 0, \quad \mbox{and} 
\\
\labell{c1Ia}
\left(\frac{2}{mn} - \frac{a}{m^2} - \frac{b}{n^2} \right)
+ \left(\frac{1}{m^2} - \frac{2}{lm} \right) +
\left( \frac{1}{n^2} + \frac{2}{ln} \right) = 0 .
\end{gather}

Assume first that  $l \geq m$ and $l \geq n$.  
If $l = 1$, then $m = n = 1$ as well.
Otherwise,  applying Lemma~\ref{modulo} to the isotropy sphere  $N \subset M^{\Z_l}$ 
implies that $m = n \mod l$, and hence again $m = n$.
Equations \eqref{1Ia} and \eqref{c1Ia} now simplify to
$a + b = 0$ and $a + b = 4$, respectively.
This is impossible.

Assume next that  $m > l$ and $m \geq n$.  
As we saw in the third paragraph of this proof,
the fact that  $m \neq 1$ implies that  $a = 1$.
Moreover,  applying Lemma~\ref{modulo} to the 
isotropy submanifold $N \subset M^{\Z_m}$
implies that  $n = l \mod m$, and hence  $n = l$.
Equations \eqref{1Ia} and \eqref{c1Ia} now simplify to
$bn + m =  0$ and $b = 3$, respectively.
Since $n$ and $m$ are both positive, this is impossible.

Finally, assume that $n > m $ and $ n > l$.
As we saw in the third paragraph of this proof,
the fact that  $n \neq 1$ implies that $b = 1$.
Moreover, applying Lemma~\ref{modulo} to  the 
isotropy submanifold $N \subset M^{\Z_n}$
implies that  $m + l = 0  \mod n$, and hence 
$n = m + l$.
Equation \eqref{1Ia} now simplifies to $a = 1$. 
In sum,  $a = b = 1$ and $n = m + l$.
Comparing with Example~\ref{proj},
statement (A) is true  for the action
$$\lambda \cdot [x_0,x_1,x_2,x_3] =
[x_0,x_1,\lambda^m x_2, \lambda^{n} x_3].$$

\noindent
{\bf Case II: The fixed set is discrete except one non-extremal fixed surface.}

Since  $\Phi$ is a perfect Morse-Bott function,
the fixed set consists of a point $p_0$ of index $0$,
a surface $\Sigma_1$ of index $2$, and a point $p_3$ of index $6$.
Since $\Sigma_1$ has index $2$,
the normal bundle to $\Sigma_1$ is the direct sum of two line
bundles with equivariant  Euler classes $mt + a u$ and $ -n t + b u$, where
$m$ and $n$ are natural numbers, $a$ and $b$ are integers,  and  $u$ is the positive
generator of $H^2(\Sigma_1;\Z)$.
After possibly replacing $\Phi$ by $-\Phi$, we may assume
that $m \geq n$.

Assume first that $m = n  = 1$.
Let $\{l_1,l_2, l_3\}$ be the weights at $p_0$,
where $l_1 \geq l_2  \geq l_3$.
Since $m = n = 1$, if there exist any isotropy submanifolds,
then each one  has minimum $p_0$ and maximum $p_3$.
Hence,
the weights at $p_3$ 
are $\{-l_1,-l_2, -l_3\}$.

Let  $\alpha = c^{S^1}_1(M) - \Gamma_{p_0} t  \in H^{2}_{S^1}(M;\Z)$.
Since $\alpha|_{p_0} = 0$,
Proposition~\ref{basis} implies that
$\alpha|_{\Sigma_1}$ is a multiple of $-t + bu$,
which is  the equivariant Euler class
of the negative normal bundle at $\Sigma_1$.
On the other hand, 
\begin{gather*}
c^{S^1}_1(M)|_{\Sigma_1} =
c_1(\Sigma_1) + (t + a u) + (-t + b u) = c_1(\Sigma_1) + (a + b)u, \quad \mbox{and} \\
\alpha|_{\Sigma_1} = 
c_1(\Sigma_1) +  (a + b) u -  \Gamma_{p_0} t = 
c_1(\Sigma_1) +  (a + b - \Gamma_{p_0} b) u + 
\Gamma_{p_0}(-t  + bu). 
\end{gather*}
Therefore, $c_1(\Sigma_1) =  - (a + b - \Gamma_{p_0}b ) u.$ 
Applying the same argument to the moment map $-\Phi$ we see that
$c_1(\Sigma_1) = - (a + b + \Gamma_{p_3} a )u$. 
Since $\Gamma_{p_3} = - \Gamma_{p_0} \neq 0$, this implies that
$a = b$.

Since $\int_M c^{S^1}_1(M) = 0$, 
Theorem~\ref{ABBV} and Remark~\ref{rmk:ABBV} imply that
$$\frac{\Gamma_{p_0}}{\Lambda_{p_0}} - 
\Gamma_{p_0}  a  + \frac{\Gamma_{p_3}}{\Lambda_{p_3}} = 0 .$$
Since $\Lambda_{p_0} = - \Lambda_{p_3}$ and  
$\Gamma_{p_0} = - \Gamma_{p_3} \neq 0$,
this simplifies to
$a = \frac{2}{\Lambda_{p_0}}.$
Hence,  either 
$ a = 1$, $l_1 = 2$, and $l_2 = l_3 = 1$; or
$a = 2$ and  $l_1 = l_2 = l_3 = 1$.
In either case, $c_1(\Sigma_1) = 2$ and so
$\Sigma_1$ is a sphere.

Suppose first that 
$a = b = 1$, and $l_1 = 2$, and $l_2 = l_3 = n =  m = 1$.
Comparing with Example~\ref{proj},
statement (A) is true  for the circle action
$$\lambda \cdot [x_0,x_1,x_2,x_3] = 
[x_0,\lambda x_1, \lambda x_2, \lambda^2 x_3].$$

Suppose instead that
$a =b = 2$ and  $l_1 = l_2 = l_3 = n = m =  1$.
Comparing with Example~\ref{gras},
statement (B)  is true for the circle action on $\Tilde{G}_2(\R^5)$
induced by the action on $\R^5 = \R \times \C^2$ given by
$$\lambda \cdot (t,x_1,x_2)   = (t,  x_1, \lambda x_2).$$

So instead,  assume that $m \neq 1$.
Then there is a $4$-dimensional isotropy submanifold
$N \subset M^{\Z_m}$  with minimum $\Sigma_1$ and maximum $p_3$.
Since $\Phi|_N$ is a perfect Morse-Bott function, Poincar\'e duality
implies that $\Sigma_1$ is  sphere.
Moreover, Lemma~\ref{fourdimb} implies that  $a = 1$.
Similarly, if $n \neq 1$, then 
there is a $4$-dimensional isotropy submanifold
$N \subset M^{\Z_n}$  with maximum $\Sigma_1$ and minimum $p_0$.
If there exists  other isotropy submanifolds, then each one has minimum $p_0$ and maximum $p_3$.
Hence, there exists a natural number $l$
so that the weights at $p_0$ and $p_3$
are  $\{n,n,l\}$ and  $\{-m,-m,-l\}$, respectively.

Since $\int_M 1 = 0$ and $\int_M c^{S^1}_1(M) = 0$, Theorem~\ref{ABBV} and Remark~\ref{rmk:ABBV} imply
 that
\begin{gather}
\labell{3IIc}
\frac{1}{n^2l} + \left(\frac{1}{m^2 n} - \frac{b}{m n^2} \right) 
- \frac{1}{m^2l} = 0, \quad \mbox{and}
\\
\labell{3IId}
\frac{2n + l}{n^2 l}   - 
\left( \frac{1}{m^2} + \frac{b}{n^2}   + \frac{2}{mn} \right) +
\frac{2m +l}{m^2 l} = 0.
\end{gather}
%Finally, since $c_1(M)^2 - 2c_2(M)$ has degree $4 < 6$,  
%$\int_M( c_1(M)^2 -2 c_2(M)= 0$.
%Therefore, Theorem~\ref{ABBV} implies that
%\equation\labell{3IIe}
% \frac{2n^2 + l^2}{n^2l }  
%+\left(
%\frac{b}{m} + \frac{n}{m^2} - \frac{bm}{n^2} - \frac{1}{n}
%\right)
%-  \frac{2m^2 + l^2}{m^2 l} = 0.
%\end{equation}

Suppose first that $m \geq  l$.
Then  applying Lemma~\ref{modulo} to the $4$-dimensional isotropy 
submanifold  $N \subset M^{\Z_m}$ implies that $n = l \mod m$, and so $n = l$.
Then equations \eqref{3IIc} and \eqref{3IId}  simplify to
$nb = m$ and $3n  = b$.
Hence, $3 n^2 = m$; but this contradicts $n \geq m$.

So assume instead that $l > m$.
Then applying Lemma~\ref{modulo} to  the isotropy
sphere $Z \subset M^{\Z_l}$ implies that $m + n = 0 \mod l$, and so  $l = m + n$.
Then equation \eqref{3IIc} simplifies to  $b = 1$.
In sum, $a = b = 1$ and the weights at $p_0$ and $p_3$ are 
$\{n,n ,  m+ n \}$ and $\{-m,-m, - (m + n)\}$.
Comparing with Example~\ref{proj},
statement (A) is true for the circle action on $\CP^3$ given by
$$ \lambda \cdot [x_0,x_1,x_2,x_3] =
[ x_0, \lambda^n x_1,\lambda^n x_2,\lambda^{m+n} x_3].
$$

\noindent
{\bf Case III: The fixed set contains more than one 
fixed surface.}

Since $\Phi$ is a perfect Morse-Bott function,
the fixed set consists of two  surfaces $\Sigma_0$ and $\Sigma_2$
of the same genus with index $0$ and $4$, respectively.
Let $u$ and $u'$ be the positive generators of $H^2(\Sigma_0;\Z)$
and $H^2(\Sigma_2;\Z)$, respectively.
Our first claim is that $\Sigma_0$ and $\Sigma_2$ are spheres.
Suppose on the contrary that $\Sigma_0$ and $\Sigma_2$ have  positive genus.
Then there exist classes $v_1$ and $v_2 \in H^1(\Sigma_0;\Z)$ so
that $v_1  v_2 =  u$.
By Proposition~\ref{basis}, 
there exist classes $\alpha_1$ and $\alpha_2
\in H^1_{S^1}(M;\Z)$ so that $\alpha_i|_{\Sigma_0} = v_i$;
let $\beta = \alpha_1 \alpha_2$.
On the one hand, 
$\beta|_{\Sigma_0} = u$,
and so  Proposition~\ref{basis} implies that $c^{S^1}_1(M) \in H^2_{S^1}(M;\Z)$
is a linear combination of $t$ and $\beta$.
On the other hand, 
because $\alpha_i|_p = 0$ for dimensional reasons,
$\beta|_p = 0$ for all fixed points $p$. 
Since 
$\Gamma_{\Sigma_0} > 0 > \Gamma_{\Sigma_2}$,
this contradicts the facts that
$c^{S^1}_1(M)|_p = \Gamma_F \, t$ 
for all $p$ in a fixed component $F$. 
Hence, $\Sigma_0$ and $\Sigma_2$ are spheres. 

The normal bundle over $\Sigma_0$ splits as the sum
of two line bundles with equivariant
Euler classes  $mt + au$ and $nt + bu$, 
where $m$ and $n$ are natural numbers, and $a$ and $b$ are integers. 
Since  every 
isotropy submanifold must have minimum
$\Sigma_0$ and maximum $\Sigma_2$, 
the negative weights at $\Sigma_2$ are $-m$ and $-n$.
Hence, the normal bundle over $\Sigma_2$ splits as the sum
of two line bundles with equivariant
Euler class  $-mt + cu'$ and $-nt + du'$,  where $c$ and $d$
are integers.

Since the degree of 
 $(c^{S^1}_1)^2 - 2 c^{S^1}_2$  is $4 < 6$,
$\int_M 1  = 0$, $\int_M c^{S^1}_1 = 0$,
and $\int_M ( (c^{S^1}_1)^2 - 2 c^{S^1}_2) = 0$.
Therefore,
Theorem~\ref{ABBV} and Remark~\ref{rmk:ABBV} imply  that
\begin{gather}
\labell{3IIIa}
 -\left ( \frac{a}{m^2n} +\frac{b}{mn^2} \right) +
\left(\frac{c}{m^2 n} + 
\frac{d}{mn^2} \right) = 0,
\\
\labell{3IIIb}
 \left(\frac{2}{mn} - \frac{a}{m^2} - \frac{b}{n^2} \right) +
 \left(\frac{2}{mn} - \frac{c}{m^2} - \frac{d}{n^2} \right)  = 0, \quad \mbox{and}
\\
\labell{3IIIc}
 \left(\frac{a}{n}   + \frac{b}{m} - \frac{na}{m^2} - \frac{mb}{n^2}
\right)
 - \left(\frac{c}{n}   + \frac{d}{m} - \frac{nc}{m^2} - \frac{md}{n^2}
\right) = 0.
\end{gather}

We may assume that $m \geq n$. If $m \neq 1$,
then applying Lemma~\ref{modulo} to  the isotropy submanifold 
$N \subset M^{\Z_m}$
implies that that $2n = 0 \mod m$.
Since $m$ and $n$ are relatively prime, this implies
that $m =2$ and $n = 1$.
Hence $m \leq 2$ and $n = 1$.

Assume first that $m = n = 1$.
Then equations \eqref{3IIIa} and \eqref{3IIIb}  simplify to
$a + b = c + d$ and $a + b + c + d = 4$. 
In sum, $m = n = 1$ and $a + b  = c + d = 2$.
 Comparing with Example~\ref{proj}, statement (A) is true for the circle action 
on $\CP^3$  given by
$$\lambda \cdot [x_0,x_1,x_2,x_3] =
[ x_0, x_1, \lambda x_2, \lambda x_3].$$

So assume instead that  $m =2$ and $n = 1$.
Applying Lemma~\ref{fourdimb} to the isotropy submanifold 
$N \subset M^{\Z_2}$, we see that $a = - c$.
Therefore,  equations \eqref{3IIIa}, \eqref{3IIIb},
and \eqref{3IIIc} simplify to  
$a + b = d$, $b + d = 2$, and $a + d = b$, respectively.
In sum, $m = 2$, $n = 1$,  $a = c = 0$ and $b = d = 1$.
Comparing with Example~\ref{gras},
statement (B) is true for the circle  action
on $\Tilde{G}_2(\R^5)$
induced by the action on $\R^5 = \R \times \C^2$
given by 
$$\lambda \cdot (t,x_1,x_2) =  (t,\lambda x_1, \lambda x_2).$$

\noindent
{\bf Case IV: The fixed set contains a 
four dimensional component.}

Let $F$ be a  four dimensional component of the fixed set.
Since $F$ is symplectic, $H^2(F;\R) \neq 0$.
Since $\Phi$ is a perfect Morse-Bott function,
this implies that $M^{S^1}$ consists of $F$ and one
isolated fixed point.
%Since these are the only fixed components, their
%are no isotropy submanifolds.
%Hence, the weights at $p_3$ are $\{-1,-1,-1\}$.
%By the local normal form theorem, there is a neighborhood of
%$p_3$ which is equivariantly symplectomorphic to the
%circle action  on $\C^3$ given by 
%$$\lambda \cdot (x_0 ,x_1,x_2) = 
%(\lambda^{-1} x_0, \lambda^{-1} x_1, \lambda^{-1} x_2).$$
%Therefore, the reduced space just below 
%$p_3$ is equivariantly diffeomorphic to
%the diagonal action of $S^1$ on $S^5$.
%regular moment fibers  are all diffeomorphic.
%But then the reduced space just above $F_0$ is a circle bundle
%over $F_0$.  Since these are equal, $F_0$ must
%be diffeomorphic to $\CP^2$, and the line bundle
%is the generator.
%
%Hence, the fixed set consists of a point $p_0$ and
%a $\CP^2$ $F_0$.  The weights at $p_0$ are $\{-1,-1,-1\}$.
%The normal bundle to $\CP^2$ is a line bundle with first
%Chern class $u + t$.
By Delzant \cite{Del} -- after possibly rescaling $\omega$ and reversing
the circle action -- $(M,\omega)$ is equivariantly
symplectomorphic to $\CP^3$ with its standard symplectic form
and the circle action given by
$$ \lambda \cdot [x_0,x_1,x_2,x_3] = [x_0,x_1,x_2,\lambda x_3].$$
A fortiori, statement (A) is true.

This completes the proof of Proposition~\ref{prop:notdiscrete}.

\section{The case that the fixed set is discrete: defining the multigraph}\labell{s:graph}

In the next remainder of the paper, we prove Theorem~\ref{thm:2} in
the case that the fixed set is discrete.  In particular,
in this section, we assume that
the fixed set is discrete
and define an associated  multigraph which is {\bf labeled}:
a real number and an even integer are
associated to each vertex, 
and a natural number $l_e$ (the {\bf length} of $e$) is  associated to each edge $e$.

Let the circle act faithfully on a $6$-dimensional compact 
symplectic manifold $(M,\omega)$ with moment map $\Phi \colon M \to \R$;
assume that the fixed  set $M^{S^1}$ is discrete. 
We define a associated labeled multigraph as follows:
The vertex set is the fixed set $M^{S^1}$;
each fixed point $p$ is labeled by it moment image $\Phi(p)$ and its index $2 \lambda_p$.
Given distinct $p$ and $q \in M^{S^1}$, let $E_{pq}$ denote the set of edges joining $p$ and $q$.
If  $\Phi(p) \leq  \Phi(q)$,  
then there is a (unique) edge $e \in E_{pq}$ of length $k \neq 1$ exactly if 
the following are all true:
\begin{enumerate}
\item $p$ and $q$ lie in the same component $N \subset M^{\Z_k}$.
\item $k$ is one of the weights at $p$; $-k$ is one of the weights at $q$.
\item  The index of $\Phi|_N$ at $q$ is equal to $2$ plus the index of $\Phi|_N$ at $p$.
\end{enumerate}
(In particular, $\Phi(p) < \Phi(q)$.)
We say the edge $e \in E_{pq}$ has {\bf minimum} $p$ and  {\bf maximum} $q$.

A multigraph contains  {\bf multiple edges} 
if there are several edges with the same minimum and maximum.
The multigraph is {\bf simple} if there are no multiple edges. 
%
%one of the following is true:
%\begin{enumerate}
%\item There is
%an isotropy sphere $Z \subset M^{\Z_k}$
%which contains both  $p$ and $q$.
%\item There is an $4$-dimensional isotropy submanifold
%$Z \subset M^{\Z_k}$ containing $p$ and $q$.  
%Moreover, $p$ is the minimum of $Z$, $q$ is not
%the maximum of $Z$,  and $-k$ is one of the weights at $q$.
%\item There is an four-dimensional isotropy submanifold
%$Z \subset M^{\Z_k}$ containing $p$ and $q$.  
%Moreover, $q$ is the maximum of $Z$, $p$ is not
%the minimum of $Z$,  and $k$ is one of the weights at $p$.
%\end{enumerate}

\begin{Example}\labell{ex:graph}
Fix natural numbers $m$, $n$, and $k$, and
consider the circle action on  $\CP^3$  given by
$$ \lambda \cdot [x_0,x_1,x_2,x_3] = [ x_0, \lambda^m x_1, \lambda^{m+n}
x_2, \lambda^{m+n+k} x_3].$$
There is a unique fixed point $p_i$ of index $2i$ for all $i$
such that $0 \leq 2i \leq 6$.
If $1 \not\in \{m,n,k\}$,
the associated multigraph  is the complete graph on  $\{p_0,p_1,p_2,p_3\}$; 
moreover,
$l_{01} = m,\  l_{12} = n,\ l_{23} = k,\  l_{02} = m + n,\ l_{13} = n + k,
\ $and $l_{03} = m + n + k,$
where $l_{ij}$ is the length of the edge from $p_i$ to $p_j$.
(In contrast if, for example, $m = 1$ then there is no
edge joining $p_0$ to $p_1$.)
\end{Example}

This multigraph has a number of nice properties.
First, the labeled multigraph determines the weights at every fixed point.
Second, the weights at the minimum and maximum of any edge  
are equal modulo its length.
Third, if two edges have the same minimum and maximum, then their lengths
are relatively prime.
To prove the first property, we need the  following proposition to
analyze the isotropy submanifolds which might arise; see \cite{Kar}.

\begin{Proposition}[Karshon]\labell{4dima}
Let the circle act faithfully on a $4$-dimensional compact
connected symplectic manifold $(N,\sigma)$ with
moment map $\Psi \colon N \to \R$; assume that the
fixed set $N^{S^1}$  is discrete.
The multiplicity of the weight  $+1$ at the minimum
is the number of fixed points  of index $2$ with negative weight $-1$.
Similarly,
the multiplicity of the weight  $-1$ at the maximum
is the number of fixed points  of index $2$ with positive weight $+1$.
\end{Proposition}

We are now ready to prove the properties described above.

\begin{Lemma}\labell{weights}
Let the circle act faithfully on a $6$-dimensional compact 
symplectic manifold $(M,\omega)$ with moment map $\Phi \colon M \to \R$.
Assume that the fixed  set $M^{S^1}$ is discrete  and consider the associated
labeled multigraph.
Given any fixed point $p$,
there are at most $\lambda_p$ edges with maximum $p$
and at most $3 - \lambda_p$ edges with minimum $p$.
Moreover, the multiset of weights at $p$
is the multiset of $\lambda_p$ negative integers and $3 -\lambda_p$
positive integers obtained by adding  $1$
and $-1$ with appropriate multiplicity to  
$$ \{ \sign(\Phi(q) - \Phi(p))\  l_e \mid  q \in M^{S^1} \mbox{ and } 
 e \in E_{pq}  \} .$$
\end{Lemma}

\begin{proof}
It is enough to show that for all $k > 1$,
the multiplicity of the weight $k$ at $p$ is the number 
of edges $e \in E$ with minimum $p$ and length $k$
(the analogous claim then holds for all $k < -1$).

Let $N$ be the component of $M^{\Z_k}$ which contains $p$.
If $\dim (N) = 0$  or $\dim (N) = 2$, 
then the claim above is obvious.
On the other hand, if   $\dim (N) = 4$  
the claim follows immediately from Proposition~\ref{4dima},
where we consider the faithful $S^1/\Z_k$ action on $N$.
\end{proof}

\begin{Lemma}\labell{compatible}
Let the circle act faithfully on a $6$-dimensional compact 
symplectic manifold $(M,\omega)$ with moment map $\Phi \colon M \to \R$.
Assume that the fixed  set $M^{S^1}$ is discrete  and consider the associated
labeled multigraph.
Given any edge $e \in E_{pq}$, 
the weights at $p$ and $q$ are equal modulo $l_e$.
\end{Lemma}

\begin{proof}
Since $p$ and $q$
are contained in the same  component of $M^{\Z_{l_e}}$, 
the claim follows immediately from Lemma~\ref{modulo}.
\end{proof}

\begin{Lemma} \labell{pairprime}
Let the circle act faithfully on a $6$-dimensional compact 
symplectic manifold $(M,\omega)$ with moment map $\Phi \colon M \to \R$.
Assume that the fixed  set $M^{S^1}$ is discrete  and consider the associated
labeled multigraph.
Given any distinct edges $e$ and $e'$ in $E_{pq}$,
the lengths $l_e$ and $l_{e'}$ are relatively prime.
\end{Lemma}

\begin{proof}
Assume not;
let $k \neq 1$ be the greatest common divisor of $l_e$ and $l_{e'}$.
We may also assume that $\Phi(p) \leq \Phi(q)$ and that $l_e \geq l_{e'}$.
Let $N \subset M^{\Z_{l_e}}$ and $N' \subset M^{\Z_{l_{e'}}}$ be
the components which contain $p$ and $q$.
The component $Z \subset M^{\Z_k}$ which contains $p$
also contains $N$ and $N'$, and hence  $q$.
By construction, the weights  at $p$ and $q$ are $\{l_e,l_{e'},a\}$ and 
$\{-l_{e},-l_{e'},b\}$, respectively,  for some integers $a$ and $b$.
Since the action is faithful, $a$ is not a multiple of $k$;
therefore,  $l_{e'} \neq a  \mod l_{e}$. 
On the other hand, by Lemma~\ref{modulo}
the  weights at $p$ are equal to the
weights at $q$ modulo $l_e$.  
Since $l_e = - l_e \mod l_e$, this implies that $2 l_{e'} = 0 \mod l_e$.
Since $l_e \geq l_{e'}$, this implies that $l_e$ is a multiple of
$l_{e'}$, that is, $l_{e'} = k$.
But then the index of $\Phi|_{N'}$ at $p$ is $0$, while the
index of $\Phi|_{N'}$ at $q$ is $4$. This is impossible
by  assumption (3) in the definition of the multigraph.
\end{proof}

\begin{Remark}\labell{extend}
We will sometimes add edges of length $1$ to $G$
in order to  reduce the number of cases that we need to consider.
(For example, this allows us to drop the condition
$1 \not\in \{m,n,k\}$ in Example~\ref{ex:graph}.)
Lemmas \ref{compatible} and \ref{pairprime} clearly 
still hold for this ``extended'' multigraph.
Moreover, as long as we add these
edges so that there are still at most $\lambda_p$ edges
with maximum $p$ and $3 - \lambda_p$ edges with minimum $p$,
Lemma~\ref{weights} still holds.
\end{Remark}

\section{The case that the fixed set is discrete and the associated multigraph is simple }
\labell{s:simple}

In this section, we prove Theorem~\ref{thm:2} in the 
case that the fixed set is discrete and the associated multigraph
is  simple, that is,  contains no multiple edges.
In fact, in this case 
only the first two possibilities can arise. 

\begin{Proposition}\labell{prop:simple}
Let the circle  act faithfully on a $6$-dimensional
compact symplectic manifold $(M,\omega)$ with moment map 
$\Phi \colon M \to \R$. 
Assume that $H^2(M;\R) = \R$,
the fixed set  is discrete,
and  the associated multigraph is simple.  
Then one of the  following two statements is true:
\begin{enumerate}
\item [(A)]
There is a subgroup $S^1 \subset \SU(4)$
and  bijection
$f \colon M^{S^1} \to \bigl(\CP^3 \bigr)^{S^1}$
so that
$T (M)|_{M^{S^1}} \cong f^*\left(
T \bigl(\CP^3 \bigr) \big|_{ ( \CP^3 )^{S^1} }
\right). $ 
\item [(B)]
There is a subgroup $S^1 \subset \SO(5)$
and  bijection
$f \colon M^{S^1} \to \Tilde{G}_2(\R^5)^{S^1} $  
so that 
$T(M)|_{M^{S^1}} \cong 
f^*\left(T \bigl(\Tilde{G}_2(\R^5) \bigr) 
\big|_{\Tilde{G}_2(\R^5)^{S^1}}\right) $. 
\end{enumerate}
%Moreover, the manifold $M$ does not contain a pair
%of isotropy spheres which intersect in two points.
\end{Proposition}

\begin{Remark}
Since the associated multigraph is simple, the manifolds described
above do not contain a pair of isotropy spheres which intersect in two points.
\end{Remark}

Our proof relies heavily on the following technical lemma.

\begin{Lemma}\labell{tech:graph}
Fix a natural number $l_{ij} \geq 1$ for each pair 
$\{i,j\} \subset \{0,1,2,3\}$.
Assume that the multisets $\{ \sign(k - i) \,  l_{ik} \mid k \neq i \}$
and $\{ \sign(k - j) \, l_{jk} \mid k \neq j \}$
are equal modulo $l_{ij}$ for each such pair.
Then one of the following three statements is true:
\begin{itemize}
\item[(a)] 
$l_{02} \leq l_{01} + l_{12}$, \,
$l_{13} \leq l_{12} + l_{23}$, \, 
$l_{03} \leq l_{01} + l_{13}$, \,  and  \,
$l_{03} \leq l_{02} + l_{23}$. 
\item [(b)] $l_{01} = l_{23}$, \,
$l_{02} \leq l_{12} + l_{03}$,\,
$l_{13} \leq l_{12} + l_{03}$,\, and \,
$l_{03}  \leq l_{01} + l_{12}$.
\item [(c)] $l_{02} = l_{13}$  and $l_{03} \leq l_{12}$.
\end{itemize}
\end{Lemma}

\begin{proof}[Proof of Proposition~\ref{prop:simple}]
By Remark~\ref{rmk:betti},
$H^{2i}(M;\R) = \R$ for all $i$ such that $0 \leq 2i \leq 6$.
Hence, since $\Phi$ is a perfect Morse function,
there is exactly one fixed point
$p_i$ of index $2i$ for all $i$ such that $0 \leq 2i  \leq 6$.
By Proposition~\ref{cor:order}, $\Phi(p_i) < \Phi(p_j)$ exactly if  $i < j$.

By assumption, the associated labeled multigraph $G$ is simple.
Add edges of length $1$ so that  $G$ is a complete graph. 
Then for each vertex $p_i$ there are exactly
$i$ edges with maximum $p_i$ and $3 - i$ edges with minimum $p_i$.
(See Remark~\ref{extend}.) 
By  Lemma~\ref{weights} 
the multiset of weights at $p_i$ is
$\{ \sign(k - i) \, l_{ik} \mid k \neq i \}$,
where $l_{ik}$ denotes  the length of the edge $\{i,k\}$.
By Lemma~\ref{compatible}, 
the  weights at $p_i$ and  $p_j$ are equal modulo $l_{ij}$
for each pair $\{i,j\} \subset \{0,1,2,3\}$.
Therefore, Lemma~\ref{tech:graph}  implies that one of 
the three statements
(a), (b), or (c) is true.
Finally, since $c^{S^1}_1(M)$ has degree $2 < 6$, $\int_M c^{S^1}_1(M)  = 0$.
Hence, Theorem~\ref{ABBV} (together with Remark~\ref{rmk:ABBV}) implies that
\begin{equation}\labell{simplec1}
\frac{l_{01} + l_{02} + l_{03}}{l_{01} l_{02} l_{03} } + 
\frac{l_{01} - l_{12} - l_{13} }{l_{01} l_{12} l_{13} } + 
\frac{l_{23} -l_{02} - l_{12} }{l_{02} l_{12} l_{23} } + 
\frac{l_{03} + l_{13} + l_{23} }{l_{03} l_{13} l_{23} }  = 0.
\end{equation}

First assume that  statement (a) is true. Rewrite equation \eqref{simplec1} as
$$ 
\frac{l_{01} + l_{12} - l_{02}}{l_{01} l_{12} l_{02}} 
+
\frac{l_{12} + l_{23} - l_{13}}{l_{12}l_{23} l_{13}}
+ 
\frac{ l_{01} + l_{13} - l_{03}}{l_{01} l_{13} l_{03}}   
+ 
\frac{l_{02} + l_{23}  - l_{03}}{l_{02} l_{23}  l_{03}} = 0
.$$
Since all the lengths are positive,
this implies that  the inequalities in statement (a)
are equalities, that is, 
$l_{02} =  l_{01} + l_{12}$, \,
$l_{13} = l_{12} + l_{23}$,  and
$l_{03} =  l_{01} +   l_{12} + l_{23}$. 
Comparing with  Example~\ref{proj},
the weights at the fixed points agree with those associated to
the circle action on $\CP^3$  given by
$$\lambda \cdot [x_0,x_1,x_2,x_3] 
= [ x_0, \lambda^{l_{01}} x_1, 
\lambda^{l_{01} + l_{12}} x_2, 
\lambda^{l_{01} + l_{12} + l_{23}} x_3].$$
Hence statement (A) is true.

Now assume that statement (b) is true.
Since $l_{01} = l_{23}$, we can 
rewrite equation \eqref{simplec1}  as
$$
\frac{l_{12} + l_{03} - l_{02}}{l_{12}l_{03}l_{02}} 
+
\frac{l_{12} + l_{03} - l_{13}}{l_{12} l_{03} l_{13}} 
+
2 \,\frac{ l_{01} + l_{12} - l_{03}}{l_{01} l_{12} l_{03}} 
= 0
.$$
Together with the inequalities in statement (b), this implies that  
$l_{01} = l_{23}$,  \,
$l_{03}  =  l_{01} + l_{12}$, and
$l_{02} = l_{13} =  l_{01} + 2 l_{12}$.
Comparing with  Example~\ref{gras},
the weights at the fixed points agree with those associated to
the circle action on $\Tilde{G}_2(\R^5)$ induced  by
the action on $\R^5 = \R \times \C^2 $ given by
$$\lambda \cdot (t,x_1,x_2) = 
(t , \lambda^{l_{12}} x_1 ,\lambda^{l_{01} + l_{12}} x_2).$$
Hence statement (B) is true.

Finally,  assume that statement (c) is true.
Since $l_{02} = l_{13}$, we can 
rewrite equation \eqref{simplec1}  as
$$ \frac{2}{l_{02} l_{03}} + \frac{2}{l_{02} l_{12}} +
\frac{ l_{12} - l_{03}}{l_{01}  l_{12} l_{03}} +
\frac{ l_{12} - l_{03}}{l_{23}  l_{12} l_{03}} = 0 .$$ 
This contradicts the fact that   $l_{03} \leq l_{12}$.
\end{proof}

%\left(\frac{1}{l_{01}l_{02}} +
%\frac{1}{l_{01}l_{03}} +
%\frac{1}{l_{02}l_{03}}\right)  +
%\left( - \frac{1}{l_{01}l_{12}} -
%\frac{1}{l_{01}l_{13}}  + 
%\frac{1}{l_{12}l_{13}} \right)+ $$ 
%$$\left(\frac{1}{l_{12}l_{02}} -
%\frac{1}{l_{02}l_{23}} -
%\frac{1}{l_{12}l_{23}}\right) +
%\left(\frac{1}{l_{03}l_{13}} +
%\frac{1}{l_{03}l_{23}} +
%\frac{1}{l_{13}l_{23}} \right) = 0.  $$

We will spend the remainder of this section proving 
Lemma~\ref{tech:graph}.
Fix a natural number $l_{ij}$ for each pair $\{i,j\} \subset \{0,1,2,3\}$;
assume that they satisfy the assumptions of the lemma.
We will think of $l_{ij}$ as labeling the edge $e_{ij}$
on the complete graph $G$ on $\{p_0,p_1,p_2,p_3\}$.
We say that $e_{ij}$ is {\bf longer} than $e_{mn}$ if 
either $l_{ij} > l_{mn}$ or if $l_{ij} = l_{mn}$ and $e_{ij}$
appears before $e_{mn}$ on this list: 
$e_{01}, e_{23}, e_{12},e_{02},e_{13},e_{03}$\footnote{
We break the ties in this way to simplify the argument.
For example, the arguments are very similar
in the case that $l_{01} = l_{12}$ and the case that $l_{01} > l_{02}$.}

Fix an edge $e_{ij}$,  and let $\{m,n\} = \{1,2,3,4\} \setminus \{i,j\}$.
We say that $e_{ij}$  is {\bf regular}  if    
$\sign(m - i) \,  l_{im} = \sign(m - j) \,  l_{jm} \mod l_{ij}  
$ and $
\sign(n - i) \,  l_{in} = \sign(n - j) \,  l_{jn} \mod l_{ij}.$  
In contrast,
we say that $e_{ij}$  is {\bf goofy}\footnote{
These terms are borrowed from surfing,
where they describe which foot is in front. 
We chose them to emphasize that there are exactly two possibilities,
that one of them is more common, and that both are equally valid.}
$\sign(m - i) \, l_{im} = \sign(n -j ) \, l_{jn} \mod l_{ij}$  
and
$\sign(n - i) \, l_{in} = \sign(m - j) \,  l_{jm} \mod l_{ij}.$  
Since $l_{ij} = -l_{ij} \mod l_{ij}$,
every edge must be regular, goofy, or both.

We will need the following fact.

\begin{Lemma}\labell{techtech}
%Fix a natural number $l_{ij} \geq 1$ for each pair $\{i,j\}
%\subset \{0,1,2,3\}$.  
%Assume that the multisets $\{ \sign(k - i) \,  l_{ik} \mid k \neq i \}$
%and $\{ \sign(k - j) \, l_{jk} \mid k \neq j \}$
%are equal modulo $l_{ij}$ for each such pair.
If $l_{03} = l_{13}$ then $l_{03} \leq l_{02} + l_{23}$.
Similarly, if  $l_{03} = l_{23}$ then
$l_{03} \leq l_{01} + l_{13}$.
\end{Lemma}

\begin{proof}
Assume that $l_{03} = l_{13}$. 
If $e_{03}$ is regular then $l_{02} + l_{23} = 0 \mod l_{03}$,
whereas if $e_{03}$ is goofy then $l_{02} + l_{03} = l_{02} + l_{13} 
= 0 \mod l_{03}$, and so $l_{02} = 0 \mod l_{03}$.
Either way, since these are all natural numbers this implies
that $l_{03} \leq l_{02} + l_{23}$.
The other claim is proved similarly.
\end{proof}

Finally, notice that the labeled graph $G'$ obtained from $G$
by exchanging $p_0$ with $p_3$ and $p_1$ with $p_2$ 
also satisfies the assumption of Lemma~\ref{tech:graph}.
Moreover, if $G'$  satisfies the conclusions of the lemma, then $G$ does 
as well.
Therefore, we can replace $G$ by $G'$ at any time.
Given this symmetry, we only need to consider  eight possible cases.

\noindent
{\bf Case Ia:  ${\mathbf e_{01}}$ is the longest edge and it is regular.}

Since $e_{01}$ is the longest edge, $l_{01}$ is greater than or equal to $l_{12}$,
$l_{02}$, $l_{13}$, and $l_{03}$.
Since $e_{01}$ is  regular
$l_{03} = l_{13} \mod l_{01}$ and $l_{02} = l_{12} \mod l_{01}$. 
Since these are all natural numbers 
this implies that
$$l_{03} = l_{13} \qquad \mbox{and} \qquad l_{02} = l_{12}.$$
By Lemma~\ref{techtech}, the displayed equations imply that (a) is true.

\noindent
{\bf Case Ib:  $\mathbf{ e_{01}}$ is the longest and  is goofy.}

By an argument similar to the first paragraph above, 
$$
l_{02} = l_{13}
\qquad \mbox{and} \qquad 
l_{03} = l_{12} 
.$$ 
Hence, (c) is true.

\noindent
{\bf Case IIa: $\mathbf{ e_{12}}$ is the longest and  is regular.}

By an argument similar to the first paragraph of  case Ia,
$$l_{01} = l_{02} \qquad \mbox{and} \qquad l_{13} =  l_{23}.$$
Since $l_{01} = l_{02}$, the edge $e_{03}$ is both goofy and regular.
Thus $l_{01} + l_{13} = 0 \mod l_{03}$, and hence
$$l_{03} \leq l_{01}  + l_{13}.$$
Together, the displayed equations imply  that (a) is true.

\noindent
{\bf Case IIb: $\mathbf {e_{12}}$ is the longest and is goofy (and not regular).}

If $l_{02} = l_{13}$ then since $e_{12}$ is the
longest edge  (c) is true;
so assume that $l_{02} \neq l_{13}.$

Since $e_{12}$ is the longest edge,
$l_{12}$ is strictly greater than $l_{01}$ and $l_{23}$ 
and greater than or equal to  $l_{02}$ and $l_{13}$. 
Since $e_{12}$ is goofy,
$l_{01} + l_{23} = 0 \mod l_{12}$ and $l_{02} + l_{13} = 0 \mod l_{12}$.  
Since these are all natural numbers and $l_{02} \neq l_{13}$, this implies that
$$l_{12} = l_{01} + l_{23}  \qquad \mbox{and} \qquad l_{12} = l_{02} + l_{13} .$$

Assume first that $e_{01}$ is longer than $e_{23}$, $e_{02}$, and $e_{13}$.
Because $e_{12}$ is not regular, $l_{01} \neq l_{02}$, so this implies that 
$l_{02} < l_{01}$, $l_{23} \leq l_{01}$, and $l_{13} \leq l_{01}$.
Since also $l_{03} \leq l_{12} = l_{02} + l_{13}$, we conclude that $l_{03} < l_{01} + l_{13}$.
If $e_{01}$ is regular, then
$l_{03} = l_{13} \mod l_{01}$, and so  $l_{03} = l_{13}$.
Since $e_{12}$ is the longest edge, by Lemma~\ref{techtech}
this implies that (a) is true.
On the other hand, if  $e_{01}$ is goofy, then $l_{02} = l_{13} \mod l_{01}$, 
which  contradicts $l_{02} \neq l_{13}$.

Up to symmetry, the only remaining possibility is that 
$e_{02}$ is longer than $e_{01}$, $e_{23}$, and $e_{13}$.
This implies that  $l_{01} < l_{02}$, $l_{23} < l_{02}$, and $l_{13} \leq l_{02}$.
Since also $l_{03} \leq l_{12} = l_{01} + l_{23}$,
we conclude that $l_{03} < l_{02} +  l_{23}$.
If $e_{02}$ is regular,
then  $l_{03} = l_{23} \mod l_{02}$, 
and so $l_{03} = l_{23}$.
Since $e_{12}$ is the longest edge, by Lemma~\ref{techtech}
this implies that (a) is true.
On the other hand, if $e_{02}$ is goofy, then $l_{01} = l_{23} \mod l_{02}$,
and so  $l_{01} = l_{23}$.  Since $e_{12}$ is the longest edge,
this implies that  (b) is true.

\noindent
{\bf Case IIIa: $\mathbf{ e_{02}}$ is the longest and  is regular.}

Since $e_{02}$ is the longest edge, $l_{02}$  is strictly greater than $l_{01}$,
$l_{23}$, and $l_{12}$, and greater than or equal to $l_{03}$.
Therefore, since $e_{02}$ is regular   $$l_{02} = l_{01} + l_{12}
 \qquad \mbox{and} \qquad l_{03} = l_{23}.$$
By Lemma~\ref{techtech}, this implies that
%
%If $e_{03}$ is regular then $l_{01} + l_{13} = 0 \mod l_{03}.$
%If  $e_{03}$ is goofy then $l_{01} + l_{03} = l_{01} + l_{23}  = 0 \mod l_{03}$,
%and hence $l_{01} = 0 \mod l_{03}$.  Either way, 
$$l_{03} \leq l_{01} + l_{13}.$$
Finally, since $l_{03} = l_{23}$,   $e_{13}$  is regular and goofy.
Thus, $l_{12} + l_{23} = 0 \mod l_{13}$,  and  so 
$$l_{13} \leq l_{12} + l_{23}.$$
Together, the displayed  equations imply that  (a) is true.

%%S I am very happy with this write up 

\noindent
{\bf Case IIIb:  $\mathbf{ e_{02}}$ is the longest and  is goofy.}

Since $e_{02}$ is the longest edge, $l_{02}$  is strictly greater than $l_{01}$,
$l_{23}$, and $l_{12}$, and greater than or equal to $l_{03}$.
Therefore, since $e_{02}$ is goofy 
$$l_{01} = l_{23}  \qquad \mbox{and} \qquad l_{02} = l_{03} + l_{12}.$$

Since $l_{01} = l_{23}$, the edge $e_{13}$ is goofy, that is, 
$l_{03} + l_{12} = 0 \mod l_{13}$.  Hence,
$$l_{13} \leq l_{03} + l_{12}.$$

Assume first that  $e_{03}$ is regular. Then 
$ l_{03} + l_{12} + l_{01} = l_{02} + l_{23} =  \mod l_{03}$,
that is, $l_{01} + l_{12} = 0 \mod l_{03}$.
Hence,  $l_{03} \leq l_{01} + l_{12}$.
Together with the displayed
equations above, this implies that  (b) is true.

So assume instead that $e_{03}$ is goofy.
Then $2 l_{01} = l_{01} + l_{23}  =  0 \mod l_{03}$.
If $l_{03} \leq l_{01}$ then  (b) is true; otherwise, $l_{03} = 2 l_{01} = l_{01} + l_{23}$.
If $e_{23}$ is regular then $l_{03} + l_{12} = l_{02} =  l_{03} \mod l_{23}$,
while if $e_{23}$ is goofy then $2 l_{23} = l_{03}  = l_{12} \mod l_{23}$.
In either case, $l_{12} = 0 \mod l_{23}$  and so $l_{23} \leq l_{12}$;
hence $l_{03}  \leq l_{01} + l_{12}$. 
As before, this implies that  (b) is true.

\noindent
{\bf Case IVa: $\mathbf{ e_{03}}$ is the longest and  is regular.}

Since $l_{03}$ is strictly greater than every other $l_{ij}$
and $e_{03}$ is  regular
$$l_{03} = l_{01} + l_{13} = l_{02} + l_{23}.$$
Since $l_{03} = l_{02} + l_{23} = l_{23} \mod l_{02}$,
the edge $e_{02}$ is  regular.  Therefore,
$l_{01} + l_{12} = 0 \mod l_{02}$.
A similar argument shows 
$l_{12} + l_{23} = 0 \mod l_{13}$.
Therefore, 
$$
l_{02} \leq l_{01} + l_{12} 
\qquad \mbox{and} \qquad 
l_{13} \leq l_{12} + l_{23} 
.$$
Together,  the displayed equations  above imply that (a) is true.

\noindent
{\bf Case IVb:  $\mathbf{ e_{03}}$ is the longest and  is goofy (and not regular).}

Since $l_{03}$ is  strictly greater than
every other $l_{ij}$ and   $e_{03}$ is goofy 
$$l_{03} = l_{01} + l_{23}  = l_{02} + l_{13}.$$

Assume first that $e_{01}$ is longer than $e_{23}$, $e_{02}$, and $e_{13}$.
If $e_{01}$ is regular then $l_{01} + l_{23} = l_{03} = l_{13}  \mod l_{01}$  
which implies that  $l_{13} = l_{23}$.  This contradicts the claim that $e_{03}$ is not regular.
Therefore,   $e_{01}$ is goofy
and so  $l_{02} = l_{13} \mod l_{01}$ and $l_{03} = l_{12} \mod l_{01}$.
Since $l_{23}, l_{02}, $ and $ l_{13}$ are all less
than or equal to $l_{01}$ and  
$l_{12} < l_{03} = l_{01} + l_{23}$, 
this implies that $l_{02} = l_{13}$ and $l_{12} = l_{23}$.
Now if  $e_{02}$ is regular then 
$2 l_{02} = l_{03} = - l_{23} \mod l_{02}$, while if
$e_{02}$ is goofy then 
$2 l_{02} = l_{03} =  l_{23} \mod l_{02}$; 
either way, $l_{23} = 0 \mod l_{02}$ and so $l_{02} \leq l_{23}$.
Since $l_{02} \leq l_{01}$ and $2 l_{02} = l_{01} + l_{23}$,
this implies that  $l_{01} = l_{02}$, 
which contradicts the assumption  that
$e_{03}$ is not regular.

Up to symmetry, 
the only remaining possibility is that  $e_{02}$ is longer
than $e_{01}$, $e_{23}$, and $e_{13}$.
If $e_{02}$ is  regular
then $l_{13} + l_{02} = l_{03} = l_{23} \mod l_{02}$,
which implies that $l_{13} = l_{23}$.  This contradicts the
assumption that $e_{03}$ is not regular. 
Therefore  $e_{02}$ is goofy and so $l_{01} = l_{23} \mod l_{02}$; 
hence,  $l_{01} =  l_{23}$.
If $e_{01}$ is regular then  $2 l_{01} = l_{03} = l_{13} \mod l_{01}$,
that is, $l_{13} = 0 \mod l_{01}$. Since $l_{13} < l_{03} = 2 l_{01}$ and
$l_{01} = l_{23}$
this implies that $l_{13} = l_{23}$. 
This contradicts the assumption that $e_{03}$ is not regular.
Therefore  $e_{01}$ is goofy and so $2 l_{01} = l_{03} = l_{12}
\mod l_{01}$, that is, $l_{12} = 0 \mod l_{01}$.   
Therefore, $l_{23} = l_{01} \leq  l_{12}$. Together with the displayed
equations above, this implies that  (b) is true.

This completes the proof of Lemma~\ref{tech:graph}.

\section{The case that the fixed set is discrete but the associated multigraph is not simple}
\labell{s:notsimple}

In this section, we prove Theorem~\ref{thm:2} in the case
that the fixed set is discrete but the associated multigraph is
not simple, that is, it contains multiple edges.
In fact, in this case only the last two possibilities 
can arise.

\begin{Proposition}\labell{prop:multiple}
Let the circle  act faithfully on a $6$-dimensional compact symplectic manifold
$(M,\omega)$ with moment map $\Phi \colon M \to \R$.
Assume that  $H^2(M,\R) = \R$, the fixed set is discrete,
and the associated  multigraph  is not simple. 
Then one of the following two statements is true:
\begin{itemize}
\item [(C)]
The fixed set consists of four points; 
the weights at these points are
%%S $\{1,2,3\}$, $\{1,-1,4\}$, $\{1,-1,-4\}$, and $\{-1,-2,-3\}$.
$$\{1,2,3\}, \{1,-1,4\}, \{1,-1,-4\}, \mbox{ and } \{-1,-2,-3\}.$$
\item [(D)]
The fixed set consists of four points; 
the weights at these points are
%%S $\{1,2,3\}$, $\{1,-1,5\}$, $\{1,-1,-5\}$, and $\{-1,-2,-3\}$.
$$\{1,2,3\}, \{1,-1,5\}, \{1,-1,-5\}, \mbox{ and } \{-1,-2,-3\}.$$
\end{itemize} 
\end{Proposition}

\begin{Remark}
Every isotropy sphere $N \subset M^{\Z_k}$ must contain
exactly one fixed point with weight $-k$ and one with weight $k$.
Therefore, the manifolds described above
must contain a pair of isotropy spheres
$N \subset M^{\Z_2}$ and $N' \subset M^{\Z_3}$ 
which intersect in two points.
\end{Remark}

Given a labeled multigraph $G$ with vertices
$\{p_0,p_1,p_2,p_3\}$,
we adopt the following notational conventions:
For any  pair $\{i,j\} \subset \{0,1,2,3\}$,
$E_{ij}$ is the set of edges from $p_i$ to $p_j$.
Moreover, the edges $e_{ij}$ and $e'_{ij}$ (if they exist)
lie in $E_{ij}$ and have lengths $l_{ij}$ and $l'_{ij}$, 
respectively. 
Our proof  relies heavily on the following technical lemma.

\begin{Lemma}\labell{techmult}
Let $G$ be a labeled multigraph with vertices $\{p_0,p_1,p_2,p_3\}$.
Assume  that the following hold:
\begin{enumerate}
\item \labell{multiple} The graph is not simple and $l_e \neq 1$ 
for every edge $e$.
\item \labell{legal} The set $\cup_{j < i} \,  E_{ij}$ contains  at most $i$ edges
and the set $\cup_{j > i} \, E_{ij}$ contains at most $3-i$ edges 
for each vertex $p_i$;
let the {\bf weights} at $p_i$ be the multiset 
of $i$ negative integers and $3-i$ positive integers
obtained by adding  $1$ and $-1$ with appropriate multiplicity to
$$ \{ \sign(j - i) \, l_e \mid 
0 \leq j \leq 3 \  \mbox{and} \ e \in E_{ij} \}. $$
\item \labell{a:modulo} 
Given any edge $e_{ij}$,  the weights at $p_i$ and  
$p_j$ agree modulo $l_{ij}$.
\item \label{prime} Given any distinct pair  of edges $e_{ij}$ and $e'_{ij}$ in $E_{ij}$,  
the lengths $l_{ij}$ and $l'_{ij}$ are relatively prime.
\end{enumerate}
Then,
after possibly exchanging $p_0$ with $p_3$ and $p_1$ with $p_2$,
one of the following three statements is true.
\begin{enumerate}
\item [(x)]
$E$ is either $ \{e_{03}, e'_{03} \}$ or 
$\{e_{03}, e'_{03} ,e_{12}\}$,  
$l_{03} = 3$, and $l_{03}' = 2$.
\item[(y)]
$E = \{e_{12}, e_{12}' ,e_{01},  e_{23}, e_{03} \}$ and
$l_{12}' = l_{01} =  l_{23} =  l_{03} = 2$.
\item [(z)]
$E = \{e_{03}, e_{02} , e_{02}', e_{23}\}$, 
$l_{03} = 4$, $l_{02} = 3$, and $l_{02}' = l_{23} = 2$.
\end{enumerate}
\end{Lemma}

\begin{proof}[Proof of Proposition~\ref{prop:multiple}]
By the first paragraph of the proof of Proposition~\ref{prop:simple},
there is exactly one fixed point $p_i$ of index $2i$ for
all $i$ such that $0 \leq 2i \leq 6$
and $\Phi(p_i) < \Phi(p_j)$ exactly if $i < j$.
%Since $H^2(M;\R) = \R$ and $(M,\omega)$ is symplectic,
%Poincar\'e duality implies that
%$H^{2i}(M;\R)  = \R$ for all $i$ such that $0 \leq 2i \leq 6$.
%Hence, since $\Phi$ is a perfect Morse function,
%there is exactly one fixed point
%$p_i$ of index $2i$ for all $i$ such that $0 \leq 2i  \leq 6$.
%By Lemma~\ref{fixedorder}, $\Phi(p_i) < \Phi(p_j)$ exactly if $i < j$.

We will now check that the associated labeled multigraph $G$
satisfies  the assumptions of Lemma~\ref{techmult}.
Since $G$ is not simple,  it satisfies assumption (1).
By Lemma~\ref{weights}, assumption (2) holds and
the weights at $p_i$ is the multiset 
of $i$ negative integers and $3-i$ positive integers
obtained by adding $1$ and $-1$ with appropriate multiplicity to
$ \{ \sign(j - i) \, l_e \mid 
0 \leq j \leq 3 \  \mbox{and} \ e \in E_{ij} \}.$
Finally, assumptions (3) and (4) hold
by Lemmas~\ref{compatible} and \ref{pairprime}, respectively. 
Therefore, Lemma~\ref{techmult} implies that 
one of the statements (x), (y),  or (z) is true.

First assume that statement (x) is true. 
For some natural number  $l_{12} \geq 1$,
the weights at $p_0,p_1,p_2,$ and $p_3$
are $\{1,2,3\},$ $\{1,-1,l_{12} \} $, $\{1,-1,-l_{12}\}$, and $\{-1,-2,-3\}$,
respectively.  
By Proposition~\ref{prop:isom} (see also Example~\ref{ex:eq}), 
as a $H^*(\CP^\infty;\Z) = \Z[t]$ module $H_{S^1}^*(M;\Z)$ is generated by
$1, \alpha_1,\alpha_2,\alpha_3$, 
where
\begin{equation}\labell{alpha}
\begin{gathered}
\alpha_1|_{p_1} = t, \
\alpha_1|_{p_2} = \frac{6+ l_{12}}{6 - l_{12}}t , \
\alpha_1|_{p_3} = \frac{12}{6 - l_{12}} t, \
\\
\alpha_2|_{p_2} =   l_{12} \, t^2, \ 
\alpha_2|_{p_3} = 6 t^2,  \
 \alpha_3|_{p_3} = 6 t^3,
\ \mbox{and} \ \alpha_i|_{p_j} = 0  \ \ \forall \ \ j < i.
\end{gathered}
\end{equation}
Since $\alpha_1$ is an integral class, this implies that 
$6 - l_{12}$ divides $12$.  In particular,  $l_{12} \neq 1$
and so  there is an isotropy sphere 
$N \subset M^{\Z_{l_{12}}}$  which contains $p_1$ and $p_2$.
Let $\gamma \in H^4_{S^1}(M;\Z)$ be the  push-forward in equivariant
cohomology of the natural generator $1 \in H^0_{S^1}(N)$.
Since $p_0$ and $p_3$ are not in $N$, $\gamma|_{p_0} = \gamma|_{p_3} = 0$.
On the other hand, for $i = 1$ or $2$,
$\gamma|_{p_i}$ is the product of the weights of the normal bundle to $N$ at 
$p_i$, that is, $\gamma|_{p_1} = \gamma|_{p_2} = -t^2$.
Comparing with equation \eqref{alpha}, this implies that
$\gamma = \frac{2}{6 - l_{12}} \alpha_2 -  \alpha_1 t$.
Since $\gamma$ is an  integral linear combination of $t^2$, 
$\alpha_1 t$, and $\alpha_2$, this implies that $6 - l_{12}$ divides $2$.
Hence, $l_{12} = 4, 5, 7$ or $8$.
Finally, Lemma~\ref{comp} implies that 
$6 = \Gamma_{p_0}  > \Gamma_{p_1} = l_{12}$.  
Hence, $l_{12} = 4$ or $5$, as required.

Now assume that statement (y) is true.
Then the weights at the fixed points 
are $\{1,2,2\}, \{2,-2,l_{12}\}, \{2,-2,-l_{12}\}$
and $\{-1,-2,-2\}$ for some natural number $l_{12}$.
Since $\int_M c^{S^1}_1(M) = 0$, Theorem~\ref{ABBV} (together with Remark~\ref{rmk:ABBV}) implies that
$\frac{5}{4} - \frac{1}{4 } - \frac{1}{4  } + \frac{5}{4} = 0$,
which is not true.

Finally, assume that statement (z) holds.
Then the weights at the fixed points 
are either 
\begin{itemize}
\item $\{2,3,4\}, \{1,1,-1\}, \{2,-2,-3\}$, and $\{-1,-2,-4\}$;
or
\item $\{-2,-3,-4\}, \{-1,-1,1\}, \{-2,2,3\}$, and $\{1,2,4\}$. 
\end{itemize}
Since $\int_M 1  = 0$, Theorem~\ref{ABBV} implies that
$\frac{1}{24} -  1 + \frac{1}{12} - \frac{1}{8} = 0$, which
is not true.
\end{proof}

We will spend the remainder of this section proving 
Lemma~\ref{techmult}.
Let $G$ be any labeled multigraph satisfying the assumptions
of the lemma.

Consider any distinct pair of edges $e_{ij}$ and $e'_{ij}$  in $E_{ij}$, where $i < j$.
By assumption (\ref{prime}), the lengths $l_{ij}$ and $l'_{ij}$ are relatively prime.
By definition, the weights at $p_i$ and $p_j$ are $\{l_{ij},l'_{ij},x\}$
and $\{-l_{ij}, -l'_{ij}, y\}$, respectively, for some integers $x$ and $y$. 
By assumption (\ref{a:modulo}), these sets agree modulo $l_{ij}$ and also
agree modulo $l'_{ij}$.
Therefore, we can apply the following facts. 

\begin{Lemma}\labell{twoedges}
Let $l \geq l' \geq 1$ 
be relatively prime natural numbers;  let
$x$ and $y$ be integers.
Assume that the multisets 
$\{l,l',x\}$ and $\{-l,-l', y \}$ are equal modulo
$l$ and also equal modulo $l'$.
Then the following hold:
\begin{enumerate}
\item [(i)] If $l' \neq 1$ then
$2l' \neq 0 \mod l$, $x \neq y \mod l$, 
$x + l' = 0  \mod l$,  and $y = l' \mod l.$
\item[(ii)]  If $l \geq x > 0$ then either
$l' = 2 = l - x$ and $x$ and $y$ are  odd, 
$l' = 1 = l - x$,  or $l' = 1$ and $x = l$.
A similar claim holds if $l \geq -y > 0$.
\item [(iii)] If $l \geq  y > 0$ 
then either 
$l' =  2 = y$ and $l$ is odd,
$l'  = 1 = y$, or $l' = 1$ and $l = 2 = y$.
A similar claim holds if $l \geq -x > 0$.
\end{enumerate}
\end{Lemma}

\begin{proof}
If $2 l' = 0 \mod l$, then since $l$ and $l'$ are  relatively
prime and $l \geq l'$, $l' = 1$ and $l \leq 2$.
Otherwise, since $2l = 0 \mod l$, 
the assumption that $\{l,l',x\}$ and $\{-l,-l',y\}$ are equal modulo $l$
implies immediately that 
$2l' \neq 0 \mod l$, $x \neq y \mod l$, 
$x + l' = 0  \mod l$,  and $y = l' \mod l.$
In particular, (i) holds.
Similarly,  if $l' > 2$ then $2l \neq  0 \mod l'$ and so
$l = y \mod l'$ and $l + x  = 0 \mod l'$.

To prove (ii),  assume that $l \geq x > 0$.
If $l' = 1$ and $l \leq 2$, then $l \geq x > 0$ implies
immediately that either 
$l' = l - x$  or $x = l$.
Otherwise,   $x + l' = 0 \mod l$ by the previous paragraph,
and $l > l'$ since $l$ and $l'$ are relatively prime.
Therefore,  $l' = l - x$.
If $l' > 2$, then  $2l \neq 0 \mod l'$ and $l + x = 0 \mod l'$,
that is, $2l = 0 \mod l'$. Since these equations give a contradiction,
$l' \leq 2$, as required.
Finally, if   $l' = 2$  then $l$ is odd since $l$ and $l'$ are
relatively prime.  Since $x = l - 2$, $x$ is odd as well.
Moreover, since $2l = 0 \mod 2$, $x + y = 0 \mod 2$, and so $y$ is also odd.

To prove (iii), assume that $l \geq y > 0$.
If $l' = 1$ and $l \leq 2$, then $l \geq y > 0$ implies 
immediately that  either  $y = 1$ or  $l = 2 = y$.
Otherwise,  by the first paragraph, $y = l' \mod l$, and so $y = l'$.
If $l' > 2$, then $ l = y \mod l'$, that is, $l = 0 \mod l'$.
Since this is impossible $l' \leq 2$, as required. 
Finally, if   $l' = 2$  then $l$ is odd since $l$ and $l'$ are
relatively prime.  
\end{proof}

Finally, notice that the labeled multigraph $G'$ obtained by exchanging $p_0$
with $p_3$ and $p_1$ with $p_2$ also satisfies the assumptions of 
Lemma~\ref{techmult}
Therefore, we can replace $G$ by $G'$ at any time.
Given this symmetry, we only need to  consider four cases.

\noindent
{\bf Case I: 
$\mathbf{ E_{02}} \boldsymbol{\neq} \boldsymbol{\emptyset} $ and $\mathbf{ E_{03}}$ contains at least two edges. }

Let $l_{03} \geq l_{03}' \neq 1$ be the length of two edges $e_{03}$ and 
$e'_{03}$ from $p_0$ to $p_3$,
and $l_{02}$ be the length of an edge $e_{02}$ from $p_0$ to $p_2$. 
Let $m$ be the positive weight at $p_2$.

Note that since $l_{03}' \neq 1$, (ii) above implies that the third weight at 
$p_0$ cannot be either $l_{03}$ or $l_{03}'$.
Hence, there cannot be any other edge of length $l_{03}$
or $l_{03}'$ with minimum $p_0$.  
A similar argument proves
that there cannot be any other edge of length 
$l_{03}$ or $l_{03}'$ with maximum $p_3$.
In particular,  if $m$ is  $l_{03}$ or $l_{03}'$,
then since $m > 1$ there is an edge  from $p_2$ to $p_3$ 
with length $l_{03}$ or $l_{03}'$; this is impossible.

Suppose first that $l_{02} \geq l_{03}$ and $l_{02} \geq m$.
Then compatibility along the edge 
$e_{02}$  implies that $m$ is either $l_{03}$ or $l_{03}'$.
As we have seen above, this is impossible.

So suppose next that $m \geq l_{03}$ and $m \geq l_{02}$.
Since $m \neq 1$ there is an edge $e_{23} \in E_{23}$
of length $m$.
Compatibility along $e_{23}$
implies that $l_{02}$ is either  $l_{03}$ or $l_{03}'$.
Again, this is impossible.

Finally, suppose that $l_{03} > l_{02}$ and $l_{03} > m$.
Since $l_{03}' \neq 1$,  by (ii)  above 
$l_{03}' = 2$ and $l_{03} = l_{03}' + l_{02}$.
But then compatibility along $e_{02}$
implies that $m = l_{03}' \mod l_{02}$.
Because  $l_{03}' + l_{02}  = l_{03} > m$ and
$l_{02} \geq l_{03}'$ (since $l_{02} \neq 1$ and $l_{03}' = 2$), this implies that $m = l_{03}'$.
Once again, this impossible.

\noindent
{\bf Case II: 
$\mathbf{ E_{02}} \boldsymbol{\neq} \boldsymbol{\emptyset}$ and $\mathbf { E_{03}}$ contains exactly one edge.}

By assumptions (\ref{multiple}) and (\ref{legal}), this implies
that -- up to symmetry -- there must be exactly two edges $e_{02}$ and $e_{02}'$ from
$p_0$ to $p_2$; we may assume that  $l_{02} \geq l_{02}$.
Let $l_{03}$ denote the length of the edge $e_{03}$ from
$p_0$ to $p_3$.
Since $l_{02}' \neq 1$, by (iii)  above 
the positive weight at $p_2$
is not $1$.
Hence, the graph also contains
an edge $e_{23}$ of length $l_{23}$  from $p_2$ to $p_3$.
By assumption (\ref{legal}),  after possibly adding an edge of length $1$,
the graph consists of these edges and
an edge $e_{13}$ of length $l_{13}$ from $p_1$ to $p_3$.

Assume first that $l_{02} \geq l_{03}$ and $l_{02} \geq  l_{23}$.
Since $l_{02}'  \neq 1$, by (ii) and (iii) above
this implies that $l_{23}$ is both even and odd, which
is impossible.

Now, assume that  that  $l_{23} \geq  l_{03}$ and 
$l_{23} \geq l_{02}$.
Then compatibility along  $e_{23}$
implies that $l_{03}$ is either equal
to $l_{02}$ or $l_{02}'$. 
Since $l_{02}' \neq 1$, this  contradicts (ii) above,
just  as in the previous case.

Now,  assume that $l_{13} \geq  l_{03}$ and $l_{13} \geq  l_{23}$. 
Since the negative weights at $p_1$ is $-1$,
compatibility along $e_{13}$  implies that
$l_{03} = 1$ or $l_{23} = 1$.
This is impossible.

Finally, assume  that $l_{03} > l_{02}$,  $l_{03} > l_{23}$, and $l_{03} \geq l_{13}$.
Then compatibility along $e_{03}$
implies that 
either  $l_{03} = l_{23} + l_{02}' =  l_{13} + l_{02}$
or $l_{03} = l_{23}  + l_{02} = l_{13} + l_{02}'$.
Since $l_{02}' \neq 1$, fact (i) above implies that
$l_{03} \neq l_{23} \mod l_{02}$. 
Hence,
$$l_{03} = l_{23} + l_{02}' =  l_{13} + l_{02}.$$
But  (i) above also implies that
$l_{03} + l_{02}' = 0 \mod l_{02}$,
and $l_{23} = l_{02}' \mod l_{02}$.
Hence, $3 l_{02}' = 0 \mod l_{02}$.
Since $l_{02}$ and $l_{02}'$ are relatively prime and $l_{02} \geq l_{02}' \neq 1$,
this implies that $$l_{02} = 3 \quad \mbox{and} \quad  l_{02}' = 2.$$
Since $l_{03} = l_{02}' \mod l_{23}$, compatibility along $e_{23}$
implies that  $l_{02} = l_{13} \mod l_{23}$; hence $4 = 0 \mod l_{23}$.
Moreover,  $l_{23} = l_{02}' \mod l_{02}$, that is, $l_{23} = 2 \mod 3$.
Combined, these  imply that $$l_{23} = 2.$$
Together, the displayed equations imply that (z) is true.

\noindent
{\bf Case III: 
$\mathbf{ E_{02}} \boldsymbol{\neq} \boldsymbol{\emptyset}$ and 
$\mathbf{ E_{03}} \boldsymbol{=} \boldsymbol{\emptyset}$.}

Since $E_{02} \neq \emptyset$, assumption (\ref{legal}) implies
that $E_{12}$ cannot contain two edges.
By assumption (1), this implies that  either $E_{02}$ or $E_{13}$ must
contain two edges.
Hence, by  assumption (\ref{legal}), after possibly adding edges of length $1$,
the multigraph contains exactly six edges:
$e_{02}, e'_{02},
e_{13}, e'_{13} , e_{01}$ and $e_{23}$.
We may assume that  $l_{13} \geq l'_{13}$ 
and that  $l_{02} \geq l'_{02}$.
By assumption (1), we cannot have $l'_{02} = l'_{13} = 1$.

Assume first   that 
$l_{01} \geq l_{02}$ and $l_{01} \geq l_{13}$.
Then compatibility along $e_{01}$
implies that $l_{02} = l_{13}$ and $l_{02}' = l_{13}'$;
hence $l_{02}' = l_{13}' \neq 1$.
Since $l_ {02}' \neq 1$, (i) above
implies that $2 l_{02}' \neq 0 \mod l_{02}$ and 
$l_{02}'  = l_{23} \mod l_{02}$. 
Since $l_{13}' \neq 1$, (i) above
implies that $l_{13'} + l_{23} = 0 \mod  l_{13}$, that is, $l_{02}' + l_{23} = 0 \mod l_{02}$.
Together, these three equations give a contradiction.

Hence, by symmetry we may assume that
$l_{02} > l_{01}$, $l_{02} > l_{23}$, and $l_{02} \geq l_{13}$.
If $l_{02}' \neq 1$,  (ii) and (iii) above imply
that $l_{23}$ is both odd and even.
Since this is impossible, $l_{02}' = 1$.
Therefore, (ii) above implies that $l_{02}=   l_{01} + 1$. 
Now, since $l_{02} \geq l_{13} \geq l_{13}'$ and $l_{13}$ and $l_{13}'$ are relatively prime,
compatibility along $e_{01}$ implies that
$l_{13}' = l_{02}' = 1$, which is impossible.

\noindent
{\bf Case IV: $\mathbf{ E_{02}} \boldsymbol{=} \mathbf{ E_{13}}
\boldsymbol{=} \boldsymbol{ \emptyset}$.}

First suppose that there are three edges from
$p_0$ to $p_3$ of length
$l_{03} \geq l'_{03} \geq  l''_{03} > 1$.
Then by (ii) above $l'_{03} \leq 2$, which is impossible.
So there are at most two edges from $p_0$ to $p_3$.

Therefore, by assumption (\ref{legal}),
after possibly adding  edges of length $1$,
$G$ contains exactly six edges:
$e_{03}, e_{03}', e_{12}, e_{12}', e_{01}$, and $e_{23}$.
We may assume that $l_{03} \geq l_{03}'$ and that $l_{12}
\geq l_{12}'$. 
By symmetry, we may also assume that $l_{01} \geq l_{23}$.
By assumption (1),  we cannot have $l'_{03} = l'_{12} = 1$.

First, assume that $l_{01} \geq l_{03}$ and 
$l_{01} \geq l_{12}$.
Then compatibility along $e_{01}$
implies that $l_{03} = l_{12}$ and $l_{03}' = l_{12}'$;
hence $l_{03}' = l_{12}' \neq 1$.
Since $l_{03}' \neq 1$,  (i) above implies that
$2l'_{03} \neq 0 \mod l_{03}$ and 
$l'_{03} + l_{01} = 0 \mod l_{03}$.
Since $l_{12}' \neq 1$,  (i) above implies that
$l_{12}' = l_{01} \mod l_{12}$, that is, $l_{03}' = l_{01} \mod l_{03}$.
Together, these  three equations give a contradiction.

Next, assume that  $l_{12} \geq l_{01} \geq l_{03}$.
Since $l_{12} \geq l_{01} \geq l_{23}$,
by (iii) above either  
$l_{12}' =  l_{01}  = l_{23} = 2$ and $l_{12}$ is odd,
or $l_{12}' = 1$.
In the first case, compatibility along $e_{01}$
implies that $l_{03} = 2$ and $l_{03}' = 1$; therefore, statement (y) holds. 
In the second case, 
compatibility along $e_{01}$
implies that $l_{03}' = 1$,  which is impossible.

Next, assume  that $l_{03} > l_{01} > l_{12}$.
Then, since $l_{03} > l_{01}$, by (ii) above
$2 \geq l_{03}' = l_{03} - l_{01}$.
Since $l_{01} \neq 1$, this implies that $l_{01} \geq l_{03}'$.
So  compatibility along $e_{01}$ implies that
$l_{12} = l_{12}' = l_{03}'$.  
Since $l_{12}$ and  $l_{12}'$ are relatively prime this implies
that  $l_{12}' = l_{03}' = 1$, which is impossible.

Finally, assume that  $l_{12} \geq l_{01}$ 
and $l_{03} > l_{01}$.
Since $l_{03} > l_{01}$, by (ii) above
either $l_{03}' = l_{03} - l_{01} = 2$
and $l_{01}$ is odd, or $l_{03}' = l_{03} - l_{01} = 1$.
In the first case, since $l_{12} \geq l_{01}$ and $l_{01}$ is odd, by (iii) above $l_{12}' = l_{01} = 1$.
Since $l_{01} \geq l_{23}$ and $l_{03} = l_{01} + 2$, this implies
that  $l_{23} = 1$ and $l_{03} = 3$ as well; therefore, statement (x)
holds.
In the second case, since  
$l_{12} \geq l_{01}$ and $l_{12}' \neq 1$, by  (iii) above 
$l_{12}' =  l_{01} = 2$, and so $l_{03} = l_{01} + 1 = 3$.
But then $l_{03}$ and $l_{03}'$ are both odd,
while $l_{12}'$ is  even.
This violates compatibility along $e_{01}$.

This completes the proof of Lemma~\ref{techmult}.

\end{document}